\documentclass[3p]{elsarticle}
\usepackage{algorithm2e}
\usepackage{hyperref}
\usepackage{bm}
\usepackage{float}
\usepackage{amsmath}
\usepackage{amsfonts}
\usepackage{graphicx}
\usepackage{xcolor}
\usepackage{comment}
\usepackage{framed}
\usepackage{mathtools}
\usepackage{amssymb}
\usepackage{bm}
\usepackage[normalem]{ulem}
\usepackage{cancel}
\usepackage{caption}
\usepackage{subcaption}
\usepackage{hyperref}
\usepackage{tikz}
\usetikzlibrary{calc}

\newtheorem{remark}{Remark}
\newtheorem{definition}{Definition}

\newcommand{\kerneldim}{n}
\newcommand{\NFloating}{\widetilde{K}}
\newcommand{\energyNorm}[1]{|||#1|||}

\newcommand{\domain}{{\widehat{\Omega}}}

\newcommand{\dirichletDomain}{\Gamma^D}
\newcommand{\neumannDomain}{\Gamma^N}

\newcommand{\interactionDomain}{\Gamma_{\text{interaction}}}

\newcommand{\R}{\mathbb{R}}

\newcommand{\family}[1]{({#1})_{k=1}^{K}}
\newcommand{\ind}{\mcX}
\DeclareMathOperator{\support}{supp}
\DeclareMathOperator{\diag}{diag}
\DeclareMathOperator{\range}{range}

\newcommand{\gammab}{\boldsymbol{\gamma}}

\newcommand{\mcF}{\mathcal{F}}
\newcommand{\mcE}{\mathcal{E}}

\newcommand{\mcX}{\mathcal{X}}

\newcommand{\mcA}{\mathcal{A}}
\newcommand{\mcT}{\mathcal{T}}

\newcommand{\mbA}{\mathbb{A}}
\newcommand{\mbB}{\mathbb{B}}

\newcommand{\mbS}{\mathbb{S}}

\DeclareMathOperator*{\argmin}{arg\,min}

\def \Phib{{\boldsymbol \Phi}}
\def \db{\mathbf{d}}
\def \eb{\mathbf{e}}
\def \fb{\mathbf{f}}
\def \gb{\mathbf{g}}
\def \ub{\mathbf{u}}
\def \wb{\mathbf{w}}
\def \vb{\mathbf{v}}
\def \xb{\bm{x}}
\def \yb{\bm{y}}
\def \Bb{\mathbf{B}}
\def \Cb{\mathbf{C}}
\def \Fb{\mathbf{F}}
\def \Gb{\mathbf{G}}
\def \Ib{\mathbf{I}}
\def \Mb{\mathbf{M}}
\def \Pb{\mathbf{P}}

\def \Zb{\mathbf{Z}}
\def \Db{\mathbf{D}}

\usepackage{color}

\newcommand{\KSub}{{K}}
\newcommand{\kSub}{k}

\begin{document}

\begin{frontmatter}

  \title{A scalable domain decomposition method for FEM discretizations of nonlocal equations of integrable and fractional type}

  \author[trier]{Manuel Klar}
  \author[lanl]{Giacomo Capodaglio}
  \author[pas]{Marta D'Elia}
  \author[snl]{Christian Glusa}
  \author[uta]{Max Gunzburger}
  \author[trier]{Christian Vollmann}

\affiliation[trier]{organization={Department of Mathematics, Universit\unexpanded{\"a}t Trier},  city={Trier}, postcode={54296}, country={Germany}}
\affiliation[snl]{organization={Center for Computing Research, Sandia National Laboratories}, city={Albuquerque}, postcode={87321}, state={NM}, country={USA}}
\affiliation[lanl]{organization={Computational Physics and Methods Group, Los Alamos National Laboratory}, city={Los Alamos}, postcode={87545}, state={NM}, country={USA}}
\affiliation[pas]{organization={Pasteur Labs},  city={Brooklyn}, postcode={11205}, state={NY}, country={USA}}
\affiliation[uta]{organization={Oden Institute for Computational Science and Engineering, University of Texas}, city={Austin}, postcode={78712}, country={USA}}

  \begin{abstract}

    Nonlocal models allow for the description of phenomena which cannot be captured by classical partial differential equations.
    The availability of efficient solvers is one of the main concerns for the use of nonlocal models in real world engineering applications.
    We present a domain decomposition solver that is inspired by substructuring methods for classical local equations.
    In numerical experiments involving finite element discretizations of scalar and vectorial nonlocal equations of integrable and fractional type, we observe improvements in solution time of up to 14.6x compared to commonly used solver strategies.

  \end{abstract}



\end{frontmatter}

\section{Introduction}
\label{sec:introduction}

Nonlocal models allow for the description of phenomena which cannot be captured by partial differential equation (PDE) models.
What separates the two classes of models are that nonlocal models are derivative free and thus allow for the interaction between pairs of points separated by a non-zero or even infinite distance, whereas for PDE models, pairs of points interact only within infinitesimal neighborhoods needed to define derivatives. Perhaps the most salient difference in the solutions of the nonlocal and PDE models is that unlike the latter, nonlocal models allow for solutions having jump discontinuities.

Applications of nonlocal models include, among others, multiscale and multiphysics systems \cite{Alali2012, Askari2008, maxmultiscalea, maxmultiscaleb}, subsurface flows \cite{Benson2000,d2021analysis,Schumer2003}, corrosion and fracture propagation in solid mechanics \cite{maxcorrosiona,maxcorrosionb,maxcorrosionc,lopez2022space,silling2000reformulation}, turbulent flows \cite{jiang,suzuki2021fractional}, phase transitions \cite{BG2019vi,burkovskab}, image processing \cite{Buades2010, d2020bilevel, Gilboa2008, Gilboa2007}, and stochastic processes \cite{Burch2014, DElia2017, Du2014jump, Meerschaert2012}. Here the listed citations provide just a subsampling of the many possible additional citations. 

In this paper, the primary focus is on nonlocal diffusion settings considered in, e.g., \cite{Gunzburger2010,du12,du13,acta20} for which
the extent of interactions between pairs of points is limited to distances no larger than a given constant $\delta>0$ that we refer to as the {\em horizon.}
%
Although this setting might be construed as being too restrictive, it suffices to illustrate our central goal which is to 
%
provide a detailed, step-by-step description of a domain decomposition (DD) algorithm for the efficient determination of finite element approximations of the solution of nonlocal models.
%
In fact, such a detailed treatment for the nonlocal diffusion setting can, for the most part, be easily extended to other nonlocal models. For example, we allude to such an extension by also briefly considering the nonlocal, derivative-free peridynamics model for solid mechanics introduced in \cite{silling2000reformulation,Silling2007} which has gained widespread popularity in, among others, fracture, corrosion, and composite material settings.

DD approaches in the nonlocal setting face challenges that do not arise in the local PDE setting. For example, consider non-overlapping DD approaches in the local PDE setting wherein a finite element grid is subdivided into covering subdomains (which are often referred to as {\em substructures}) that overlap only at the common boundaries of pairs of subdomains. As a result, in this setting, {\em points located within the interior of a substructure do not interact with points in the interior of any other substructures.} As such, local subdomain problems can be defined that are independent of other subdomains except that information is exchanged between abutting subdomains only at their common boundaries.  

In our approach for DD in the nonlocal setting, we also begin by constructing the same type of non-overlapping geometric decomposition of the finite element grid. However, because of nonlocality, {\em now points in the interior of a subdomain interact with points in the interior of abutting subdomains} so that subdomain problems exchange information not just at their common boundaries, but also between pairs of points in the interiors of abutting subdomains. Thus, although the subdomains {\em geometrically} overlap only at their common boundaries, the subdomain nonlocal problems also overlap with interior points located within abutting subdomains.

Another challenge faced by DD for nonlocal problems is that, because of nonlocality, {\em the FEM stiffness matrices for both the parent single-domain and for the individual substructures are less sparse compared to those for local PDE problems.} In the local case, a given finite element interacts with only abutting fine elements whereas, in the nonlocal case, a given finite element also interacts with finite elements that are within a $\delta$-dependent neighborhood of the given element.

A third challenge faced by DD for nonlocal problems is that {\em the condition numbers of both the single-domain and multi-domain finite element stiffness matrices generally not only depend on the grid size, but also depend on the horizon $\delta$.} Due to ill-conditioning of the matrix specialized solvers and preconditioners are required to solve nonlocal equations.

In \cite{capodaglio2020general} we provided a general framework for substructuring-based DD methods for models having nonlocal interactions. Here, we develop a specific nonlocal DD methodology by combining that framework with a nonlocal generalization of the well-known FETI (Finite Element Tearing and Interconnecting) approach for defining scalable DD methods for PDEs \cite{farhat1991method,chan1994domain,mathew2008domain,toselli2006domain}. As such, here we provide the first use of the FETI approach applied to a finite element discretization of a nonlocal model. In so doing, we provide a detailed recipe for the practical implementation of a nonlocal substructuring-based method in the finite element context. Paper \cite{xu2021feti} illustrates a FETI technique for meshfree discretizations of nonlocal models. 

Section \ref{sec:singledomain} is devoted to defining both strong and weak formulations of single-domain nonlocal problems, focusing on nonlocal diffusion models and, to a lesser extent, also on nonlocal solid mechanics models. We then introduce finite element discretizations of the single-domain nonlocal problems with the aim of defining discretized weak formulations of the nonlocal problems. We provide sufficient details so that the reader becomes aware of several pitfalls faced in the process of determining approximate finite element solutions. 

Section \ref{sec:muultipledomain} is devoted to a detailed and precise discussion of general DDs for nonlocal models, again focusing on nonlocal diffusion. This process starts with same finite element grid as that used for the discretization of a single-domain nonlocal problem.
That grid is then partitioned into subdomains that overlap only at their common boundaries exactly in a similar manner as that used for overlapping DD methods in the PDE setting.
Then, for each of the overlapping substructures, discretized weak formulations of the nonlocal model are defined.
Necessarily, in particular, solutions of discretized weak formulations are defined within the overlap of two or more substructures, i.e., multiple finite element solutions are obtained within such overlaps.
Of course, globally, such solutions should be uniquely defined so that constraints have to be added that force all solutions within the overlaps to be equal to each other.

Having in hand the discretized weak formulations for multiple substructures, we apply the aforementioned extension of well-known FETI-type algorithms for PDEs to the nonlocal DD setting constructed in Section \ref{sec:muultipledomain}. As such, there results an efficient means for determining solutions of multiple-domain DD problems. 

In Section \ref{sec:FETI} a complete framework is provided for efficiently determining solutions of multiple-domain DD for nonlocal models. However, practical implementation issues related to that framework are left to Section \ref{sec:implementation}. Then, in Section \ref{sec:numerical-experiments}, computational results of some illustrative examples are provided. In particular, scalability and accuracy studies are illustrated. Finally, concluding remarks are provided Section \ref{sec:conclusion}.

\section{Nonlocal models for diffusion and solid mechanics: the single domain setting}\label{sec:singledomain}

In this section we define strong formulations and related weak formulations for two classes of nonlocal models, one being applicable to nonlocal diffusion models and the other to nonlocal solid mechanics models. We also define finite element discretizations of the weak formulations.

\subsection{Strong and weak formulations of nonlocal diffusion and mechanics models}\label{strongsingledomain}


\subsubsection{Kernel and kernel function properties}

We differentiate between nonlocal problems having scalar-valued solutions $u(\xb)\colon$ $\R^d\to\R$  and vector-valued solutions $\ub(\xb)\colon \R^d\to\R^d$.

Let $\gamma(\xb, \yb)\colon \R^d\times\R^d \to \R$ denote a compactly supported scalar-valued function, which we refer to as the {\em kernel}, having the form
\begin{equation}
\label{eq:scalarkernel}
\gamma(\xb, \yb) = \phi(\xb, \yb)\ind_{B_{\delta,p}(\xb)}(\yb),
\end{equation}
where $\ind_{B_{\delta,p}(\xb)}$ \label{inline:Bdeltap} denotes the indicator function of the $\ell^p$-ball $B_{\delta,p}(\xb)$ of radius $\delta$ centered at $\xb$. Likewise, let $\gammab(\xb, \yb)\colon \R^d\times\R^d \to \R^{d\times d}$ denote a compactly supported matrix-valued function having the form
\begin{equation}
\label{eq:vectorkernel}
\gammab(\xb, \yb) = \Phib(\xb, \yb)\ind_{B_{\delta,p}(\xb)}(\yb) .
\end{equation}
We refer to $\phi(\xb, \yb)$ and $\Phib(\xb, \yb)$ as {\em kernel functions} and we assume that they satisfy the following properties:
\begin{equation}
\label{eq:kernelfuncprop}
\begin{aligned}
&
\qquad\mbox{a. $\phi(\xb, \yb)$ and $\Phib(\xb, \yb)$ are symmetric in the function sense,}
\\[-.5ex]&
\qquad\qquad
\mbox{i.e., $\phi(\xb, \yb) = \phi(\yb, \xb)$ and $\Phib(\xb, \yb) = \Phib(\yb, \xb)$}
\\[.5ex]&
\qquad\mbox{b. $\phi(\xb, \yb)$ is non-negative and $\Phib(\xb, \yb)$ is a symmetric}
\\[-.5ex]&
\qquad\qquad\mbox{positive semi-definite matrix.}
\end{aligned}
\end{equation}
Wider classes of kernel functions are also considered in the literature; see, e.g., \cite{medu} for some examples.
A commonly found feature of kernel functions of interest is the presence of a singularity for \(\xb=\yb\).
We also note that because $\yb\in B_{\delta,p}(\xb)$ necessarily implies that $\xb\in B_{\delta,p}(\yb)$, we have that
\begin{equation}
\label{eq:symindicator}
\mbox{the indicator function $\ind_{B_{\delta,p}(\xb)}(\yb)$ is a symmetric function, i.e., $\ind_{B_{\delta,p}(\xb)}(\yb)=\ind_{B_{\delta,p}(\yb)}(\xb)$.}
\end{equation}
Combining \eqref{eq:kernelfuncprop} and \eqref{eq:symindicator}, we obviously obtain that the kernels have the following properties:
\begin{equation}
\label{eq:kernelprop}
\begin{aligned}
&
\qquad\mbox{a. $\gamma(\xb, \yb)$ and $\gammab(\xb, \yb)$ are symmetric in the function sense}
\\&
\qquad\mbox{b. $\gamma(\xb, \yb)$ is non-negative over its support and $\gammab(\xb, \yb)$ is a}
\\[-.5ex]&
\qquad\qquad\mbox{symmetric positive semi-definite matrix over its support.}
\end{aligned}
\end{equation}

\subsubsection{Geometric entities}

We are given an open, bounded subset $\domain \subset \R^d$ and are also given a bounded constant $\delta > 0$ which is referred to as the \textit{horizon}. We then define the closed {\em interaction region} $\interactionDomain$ corresponding to $\domain \subset \R^d$ by
\begin{equation}\label{interdomain}
\begin{aligned}
\interactionDomain
&
=
\big\{
\yb \in \R^d \setminus \domain
\,:\,\mbox{there exists} \;\; \xb \in \overline{\domain}\;\; \textrm{ such that }\;\; |\xb-\yb|_{p} \leq \delta
\big\}
\\&
=
\cup_{\xb \in \overline{\domain}}\Big( B_{\delta,p}(\xb)\cap (\R^d \setminus \domain) \Big),
\end{aligned}
\end{equation}
where $|\cdot|_{p}$ denotes the $\ell^p$-norm in $\R^d$, $B_{\delta,p}(\xb)\in\R^d$ denotes the $\ell^p$-ball having radius $\delta$ centered at a point $\xb$, and $\overline{(\cdot)}$ denotes the closure. Of particular interest are the values $p=2$ and $\infty$.

\begin{figure}[h!]
	\centering
	\begin{tabular}{ccc}
		\includegraphics[width=1.5in]{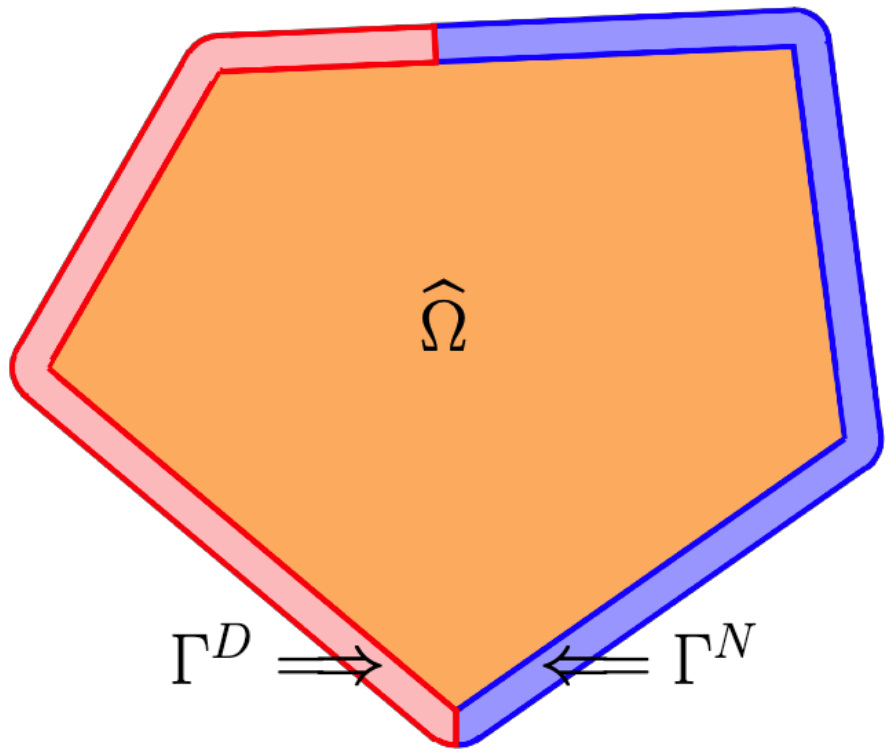}
		&
		\includegraphics[width=1.5in]{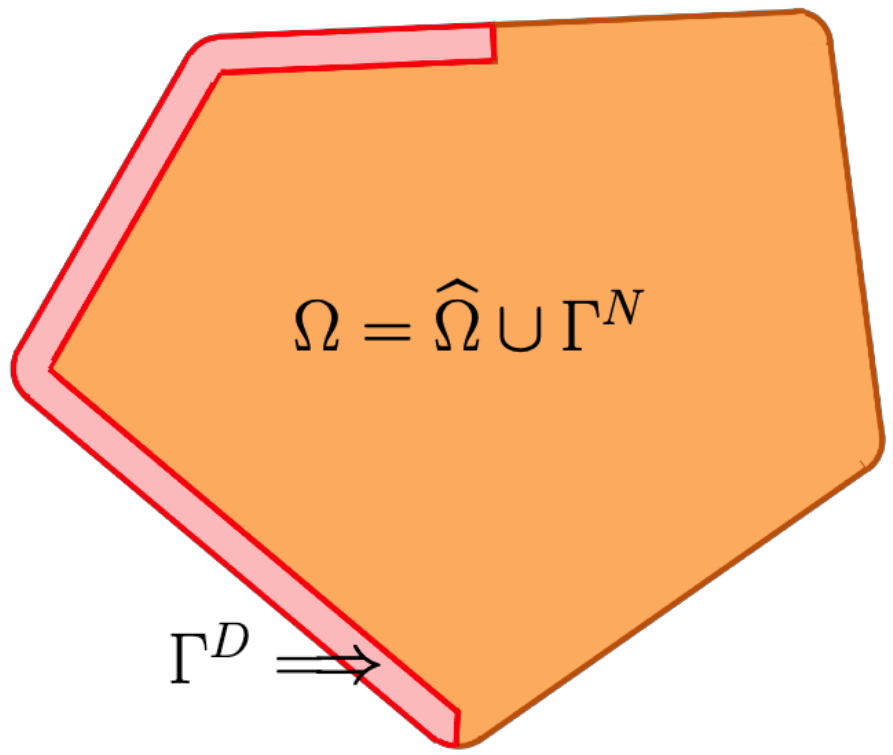}
		&
		\includegraphics[width=1.5in]{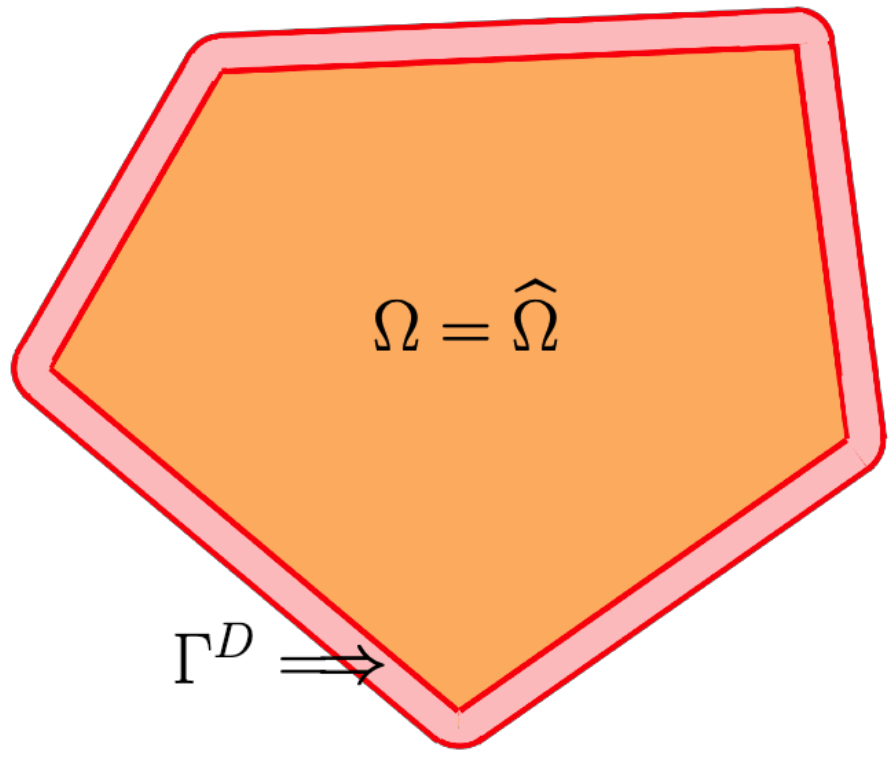}
		\\
		(a) & (b) & (c)
	\end{tabular}
	\caption{(a) An open, bounded region $\domain\subset\R^2$ and its associated interaction region $\dirichletDomain\cup\neumannDomain\subset\R^2$. (b) The region $\Omega=\domain\cup\neumannDomain$. (c) The geometric configuration in case only Dirichlet volume constraints are imposed.}
	\label{fig1}
\end{figure}

\begin{figure}[h!]
	\centering
	\includegraphics[width=2.in]{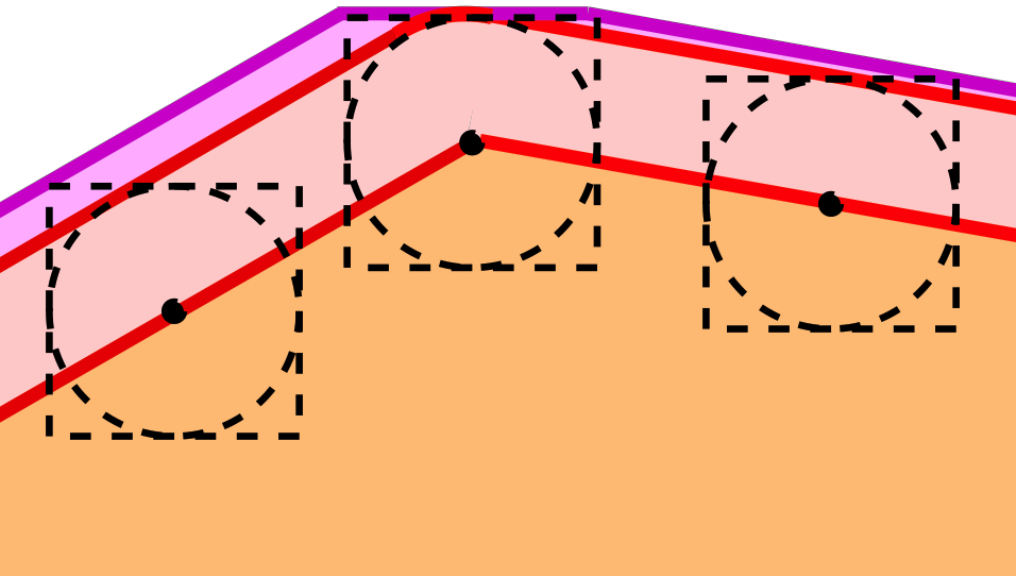}
	\caption{The orange region is a part of $\domain$ and the pink region is part of the interaction region $\dirichletDomain$ defined using $\ell^2$-balls. The pink and magenta regions together constitute part of the interaction region $\dirichletDomain$ defined using $\ell^\infty$-balls having the same radius as the  $\ell^2$-ball. Note that the width of $\dirichletDomain$ induced by the $\ell^2$-ball is constant whereas for the $\ell^\infty$-ball, that width varies with the angle the boundary of $\widehat\Omega$ makes with, say, the $x$ axis.}
	\label{fig1a}
\end{figure}

The interaction region $\interactionDomain$ is split into two disjoint parts $\dirichletDomain$ and $\neumannDomain$ \label{inline:nonlocalBoundary}. Here, $\dirichletDomain$ is a nonempty closed region whereas $\neumannDomain$ is open along its common boundary with $\dirichletDomain$ but is closed with respect to the rest of its boundary. Clearly, $\dirichletDomain\cup\neumannDomain=\interactionDomain$ and $\dirichletDomain\cap\neumannDomain=\emptyset$. Note that the width of $\dirichletDomain$ and $\neumannDomain$ depends on $\delta$ whereas $\domain$ is defined independently of $\delta$. We also define the region $\Omega=\domain\cup\neumannDomain$ \label{inline:Omega} which is open with respect to $\dirichletDomain$ but is closed with respect to the rest of its boundary. Figures \ref{fig1}a and \ref{fig1}b illustrate, for $d=2$, the four geometric entities $\domain$, $\dirichletDomain$, $\neumannDomain$, and $\Omega$. Figures \ref{fig1}c illustrates the geometric configuration in case only Dirichlet volume constraints are imposed, i.e., if $\neumannDomain=\emptyset$; obviously, in this case we have that $\Omega=\domain$.
\begin{remark}\label{rem1}
{\em
In this paper, to avoid having to deal with the approximation of curved domains $\domain$ when constructing grids, we assume that $\domain$ is a polygonal domain. However, as illustrated in Figures \ref{fig1} and \ref{fig1a}, the resulting interaction region $\interactionDomain=\dirichletDomain\cup\neumannDomain$ constructed using $\ell^2$-ball interaction regions necessarily has a curved outer boundary. In Section \ref{sec:muultipledomain} we discuss the implications that this observation has on the construction of finite element grids. On the other hand, an interaction region $\dirichletDomain\cup\neumannDomain$ constructed using $\ell^\infty$-balls does result in a polygonal outer boundary, as is illustrated in Figure \ref{fig1a}.
\qquad$\Box$}
\end{remark}
In this paper, because of their ubiquitous use, we focus on interaction regions induced by $\ell^2$- and \(\ell^{\infty}\)-balls.
As the remainder of the theoretical discussions is independent of the precise form of the kernels, we refer the reader to Section~\ref{sec:kernels} for the kernels used in numerical examples.

\subsubsection{A strong formulation of a nonlocal diffusion problem}

We are given functions $f_{\domain}(\xb)$, $f_{\neumannDomain}(\xb)$, and $g(\xb)$ defined for $\xb\in\domain$, $\xb\in\neumannDomain$, and $\xb\in\dirichletDomain$, respectively. Also given is a kernel $\gamma(\xb,\yb)$ satisfying \eqref{eq:kernelprop} for $\xb,\yb\in \domain\cup\dirichletDomain\cup\neumannDomain=\Omega\cup\dirichletDomain$. Then, a strong formulation of a nonlocal diffusion problem is given by \cite{acta20,du19,du12,du13,Gunzburger2010}.
\begin{equation}\label{eq:strong_nonloc1}
\left\{\begin{tabular}{lll}
$\displaystyle
-2\int_{\domain\cup\dirichletDomain\cup\neumannDomain}
\big(u(\yb) -u(\xb)\big)\gamma(\xb,\yb)d\yb$
$= f_{\domain}(\xb)$
& for $\xb\in\domain$,
& \qquad(a)
\\[4mm]
$u(\xb) = g(\xb)$
& for $\xb\in\dirichletDomain$,
& \qquad(b)
\\[1mm]
$2\displaystyle\int_{\domain\cup\dirichletDomain\cup\neumannDomain}
\big(u(\yb) -u(\xb)\big)\gamma(\xb,\yb)d\yb$
$= f_{\neumannDomain}(\xb)$
& for $\xb\in\neumannDomain.$
& \qquad(c)
\end{tabular}\right.
\end{equation}
The nonlocal problem \eqref{eq:strong_nonloc1} is referred to as nonlocal Poisson problem because, by analogy, (\ref{eq:strong_nonloc1}a) plays the role of a second-order elliptic partial differential equation and (\ref{eq:strong_nonloc1}b) and (\ref{eq:strong_nonloc1}c) play the roles of Dirichlet and Neumann boundary conditions, respectively. In fact, as $\delta\to0$, the problem \eqref{eq:strong_nonloc1} with a proper $\delta$-dependent scaling (see Section \ref{sec:kernels}), reduces to a boundary-value problem for a second-order elliptic partial differential equation with mixed Dirichlet and Neumann boundary conditions; see, e.g., \cite{acta20,du19,du12,du13,Gunzburger2010} for details. Clearly, (\ref{eq:strong_nonloc1}b) and (\ref{eq:strong_nonloc1}c) are not conditions that are applied on a bounding surface, but are instead applied on subsets of $\R^d\setminus\domain$ having $\delta$-dependent finite measure in $\R^d$. As such, we refer to those two conditions as {\em Dirichlet} and {\em Neumann volume constraints}, respectively.

Clearly, using $\Omega=\domain\cup\neumannDomain$, one can combine (\ref{eq:strong_nonloc1}a) and (\ref{eq:strong_nonloc1}c) into a single equation so that \eqref{eq:strong_nonloc1} can be more compactly written as
\begin{equation}\label{eq:strong_nonloc}
\left\{\begin{tabular}{lll}
$\displaystyle
-2\int_{\Omega\cup\dirichletDomain}
\big(u(\yb) -u(\xb)\big)\gamma(\xb,\yb)d\yb =f(\xb)$
\\[1.5ex]
$
\displaystyle\qquad\qquad\qquad\qquad
=f_{\domain}(\xb)\ind_{\domain}(\xb)
-f_{\neumannDomain}(\xb)\ind_{\neumannDomain}(\xb)$
& for $\xb\in\Omega=\domain\cup\neumannDomain$,
& \qquad(a)
\\[1ex]
$u(\xb) = g(\xb)$
&  for $\xb\in\dirichletDomain.$
& \qquad(b)
\end{tabular}\right.
\end{equation}
This strong formulation corresponds to the geometric configuration in Figure \ref{fig1}b and which is considered in detail in the remainder of the paper.

\begin{remark}
  If the kernel function has a singularity of high enough order (e.g. for fractional kernels) the integrals in \eqref{eq:strong_nonloc1} and \eqref{eq:strong_nonloc} can only be taken in the principle value sense.
  We do not explicitly make that distinction for the sake of a more uniform presentation.
\end{remark}

\subsubsection{A strong formulation of a nonlocal solid mechanics model}

We rename the regions $\dirichletDomain$ by $\Gamma_{\text{displacement}}$ and $\neumannDomain$ by $\Gamma_{\text{traction}}$. Then, given vector-valued functions $\fb_{\domain}(\xb)$, $\fb_{\text{traction}}(\xb)$, and $\gb_{\text{displacement}}(\xb)$ defined for $\xb\in\domain$, $\xb\in\Gamma_{\text{traction}}$, and  $\xb\in\Gamma_{\text{displacement}}$, respectively, a {\em nonlocal model for solid mechanics} that determines the vector-valued displacement $\ub(\xb)$ for $\xb\in \domain\cup \Gamma_{\text{traction}}\cup\Gamma_{\text{displacement}}$ is given by \cite{du13peri}
\begin{equation}\label{eq:strong_nonloc-v}
\left\{\begin{tabular}{lll}
$\displaystyle
-2\int_{\domain\cup\Gamma_{\text{displacement}}\cup\Gamma_{\text{traction}}}
\gammab(\xb,\yb)\big(\ub(\yb) -\ub(\xb)\big)d\yb
= \fb_{\domain}(\xb)$
&
for $\xb\in\domain$,
&
(a)
\\[3ex]
$\ub(\xb) = \gb_{\text{displacement}}(\xb)$
&
for $\xb\in\Gamma_{\text{displacement}}$,
&
(b)
\\[0ex]
$2\displaystyle\int_{\domain\cup\Gamma_{\text{displacement}}\cup\Gamma_{\text{traction}}}
\gammab(\xb,\yb)(\xb,\yb)\big(\ub(\yb) -\ub(\xb)\big)d\yb$
$= \fb_{\text{traction}}(\xb)$
&
for $\xb\in\Gamma_{\text{traction}}.$
&
(c)
\end{tabular}\right.
\end{equation}
We refer to \eqref{eq:strong_nonloc-v} as a nonlocal solid mechanics problem because, for appropriate choices of the matrix-valued kernel $\gammab(\xb,\yb)$, by analogy, (\ref{eq:strong_nonloc-v}a) plays the role of a partial differential equation model for solid mechanics, (\ref{eq:strong_nonloc-v}b) plays the role of a displacement boundary condition, and (\ref{eq:strong_nonloc-v}c) plays the role of a traction boundary condition. In fact, in such cases, as $\delta\to0$, the problem \eqref{eq:strong_nonloc-v} reduces to a boundary-value problem for a partial differential equation modeling linear elasticity with mixed displacement and traction boundary conditions; see, e.g., \cite{du13peri} for details. Again, (\ref{eq:strong_nonloc-v}b) and (\ref{eq:strong_nonloc-v}c) are referred as {\em displacement} and {\em traction volume constraints}, respectively.

In analogy with the nonlocal diffusion problem, (\ref{eq:strong_nonloc-v}a) and (\ref{eq:strong_nonloc-v}c) can be combined into a single equation.

\begin{remark}
{\em As is the case for the local partial differential equation setting, because we are dealing with vector-valued solutions $\ub(\xb)$, there are other possible types of volume constraints that can be imposed in addition to those given in \eqref{eq:strong_nonloc-v}, as is the case for the local partial differential equations setting. For example, one could impose a nonlocal Dirichlet volume constraint on the component of $\ub$ that is normal to the boundary of $\widehat\Omega$ and a nonlocal Neumann volume constraint for the tangential components. Dealing with this and other cases is beyond the scope of this paper.}\qquad$\Box$
\end{remark}

\subsection{A weak formulation of the nonlocal diffusion problem}\label{weaksingledomain}

We illustrate the construction of a weak formulation in the nonlocal diffusion setting. The construction process for the nonlocal solid mechanics setting is effected in an entirely analogous manner.

Using the nonlocal calculus developed in, e.g., \cite{Gunzburger2010} and in particular a nonlocal Green's first identity, a weak formulation of the nonlocal diffusion problem \eqref{eq:strong_nonloc} can be defined in the same manner as in case of weak formulations for partial differential equation problems for diffusion. The reader can consult these papers for details concerning how the weak formulations considered in this work are defined.

For two scalar-valued functions $u,v$, we define the bilinear form $\mathcal{A}(\cdot,\cdot)$ and linear functional $\mcF(\cdot)$ as
\begin{equation}\label{eq:bil1}
\left\{ 
\begin{aligned}
&\mathcal{A}(u,v) = 
\int_{\Omega\cup\Gamma^D}  \int_{\Omega\cup\Gamma^D}\big(u(\yb) - u(\xb)\big)\big(v(\yb) - v(\xb)\big)\gamma(\xb,\yb)d\yb  d\xb 
\\
&\mcF(v) = \int_\Omega v(\xb)f(\xb)d\xb, 
\end{aligned}\right.
\end{equation}
where $\gamma(\xb,\yb)$ and $f(\xb)$ are defined in \eqref{eq:scalarkernel} and (\ref{eq:strong_nonloc}a), respectively, and $\Omega=\widehat{\Omega}\cup\Gamma^N$ \label{inline:widehatOmega}. 
\begin{definition}
	\label{globalspace}
	\label{def:energySpace}
	We 
 define the function spaces
	\begin{align}	
	W(\Omega \cup \Gamma^D) = 
	\lbrace
	&
	w \in L^2 (\Omega \cup \Gamma^D) 
	~\colon~  \energyNorm{w}^2 < \infty\}
	\\
	W^0(\Omega, \Gamma^D) = 
	\lbrace
	&
	w \in W(\Omega \cup \Gamma^D)
	~\colon~ 
	w  = 0 
	\textrm{ in } \Gamma^D \rbrace,
	\end{align}
	where 
	\begin{align*}
	\energyNorm{w}^2 = 
	\int_{\Omega \cup \Gamma^D} 
	\int_{\Omega \cup \Gamma^D}
	(w(\yb) - w(\xb))  ^2
	\gamma(\xb, \yb) 
	d\yb d\xb 
	+ 
	\| w\|_{L^2(\Omega \cup \Gamma^D)}^2.
	\end{align*}
\end{definition}
We also define the nonlocal dual space $W'(\Omega, \Gamma^D)$ of bounded linear functionals on $W^0(\Omega, \Gamma^D)$ and, \label{inline:dualSpace}
via restriction to $\Gamma^D$, a nonlocal ``trace'' space $W^t(\Gamma^D) := \lbrace w|_{\Gamma^D} \colon w \in W(\Omega \cup \Gamma^D) \rbrace$. \label{inline:traceSpace}
Then, a weak formulation of the nonlocal volume-constrained problem \eqref{eq:strong_nonloc} is defined as follows \cite{acta20,du19,du12,du13,Gunzburger2010}:
\begin{equation}\label{eq:weak_nonloc}
\begin{aligned}
&\mbox{\it given 
	$f\in W'(\Omega, \Gamma^D)$,
	$g\in  W^t({\Gamma^D})$, and a kernel $\gamma(\xb,\yb)$,}
\\
&\mbox{\it find $u\in W(\Omega \cup \Gamma^D)$ such that}
\\
&  \qquad 
\mathcal{A}(u,v) = \mcF(v) \quad \forall\, v \in W^0(\Omega, \Gamma^D)
\\&\mbox{\it subject to $u(\xb)=g(\xb)$ for $\xb\in\Gamma^D$}. 
\end{aligned}
\end{equation}
The well-posedness of  problem \eqref{eq:weak_nonloc} is proved in, e.g., \cite{acta20,du19,du12,du13,Gunzburger2010} for $\Gamma^D$ with non-zero measure. 
\begin{remark}
{\em
The definition of specific instances of the function spaces 
$W(\Omega \cup \Gamma^D)$, 
$W'(\Omega, \Gamma^D)$, and $W^t({\Gamma^D})$
depend on integrability properties of the kernel function $\phi(\xb,\yb)$ introduced in \eqref{eq:scalarkernel}. Because the specific choices of these spaces are not germane to the goals of this paper, here we do not need to provide details about these input spaces. The interested reader may consult \cite{acta20,du19,du12,du13,Gunzburger2010} for detailed discussions about these spaces for the specific kernel functions introduced in Section \ref{sec:kernels}.
}\qquad$\Box$
\end{remark}
Problem \eqref{eq:weak_nonloc} can equivalently be expressed as an energy minimization principle \cite{acta20,du19,du12,du13,Gunzburger2010}.
To that end let 
\begin{equation}\label{singleenergy}
	\mathcal{E}(v) = \frac{1}{2}\mathcal{A}(v,v) - \mcF(v).
\end{equation}
We can then formulate \eqref{eq:weak_nonloc} as 
\begin{equation}\label{eq:minimiz_nl}
\begin{aligned}
&\mbox{\it given $f\in W'(\Omega, \Gamma^D)$, $g\in  W^t({\Gamma^D})$, and a kernel $\gamma(\xb,\yb)$,}
\\
&\mbox{\it find $u \in W(\Omega \cup \Gamma^D)$ such that}
\\
&  \qquad\qquad
\mathcal{E}(u) = \inf\limits_{v\in W} {\mathcal E}(v) \\
&\mbox{\it subject to $u(\xb)=g(\xb)$ for $\xb\in\Gamma^D$}.
\end{aligned}
\end{equation}

\subsection{Finite element discretization of the weak formulation}\label{fem-approx}

\subsubsection{Finite element grids}\label{fem-grid-onedomain}

The domain $\domain\cup\dirichletDomain\cup\neumannDomain$ is subdivided into a covering, non-overlapping set of finite elements which respects the common boundaries between the three regions by which we mean that there are no finite elements straddling across those common boundaries, a finite element vertex is located at each corner of $\domain$, and there are no hanging nodes, i.e., finite element vertices are also vertices of abutting finite elements. In so doing, we have to deal with the issue discussed in Remark \ref{rem1}, namely that even though we only consider polygonal $\domain$, the interaction region $\dirichletDomain\cup\neumannDomain$ constructed using $\ell^p$--ball interaction neighborhoods can have a curved boundary. For example, this is the case for \(\ell^{2}\)-ball interaction neighborhoods.
Similar issues arise in the evaluation of the bilinear form using quadrature and will be discussed in the next section.

\begin{figure}[h!]
\centering
\begin{tabular}{ccc}
		\includegraphics[width=1.in]{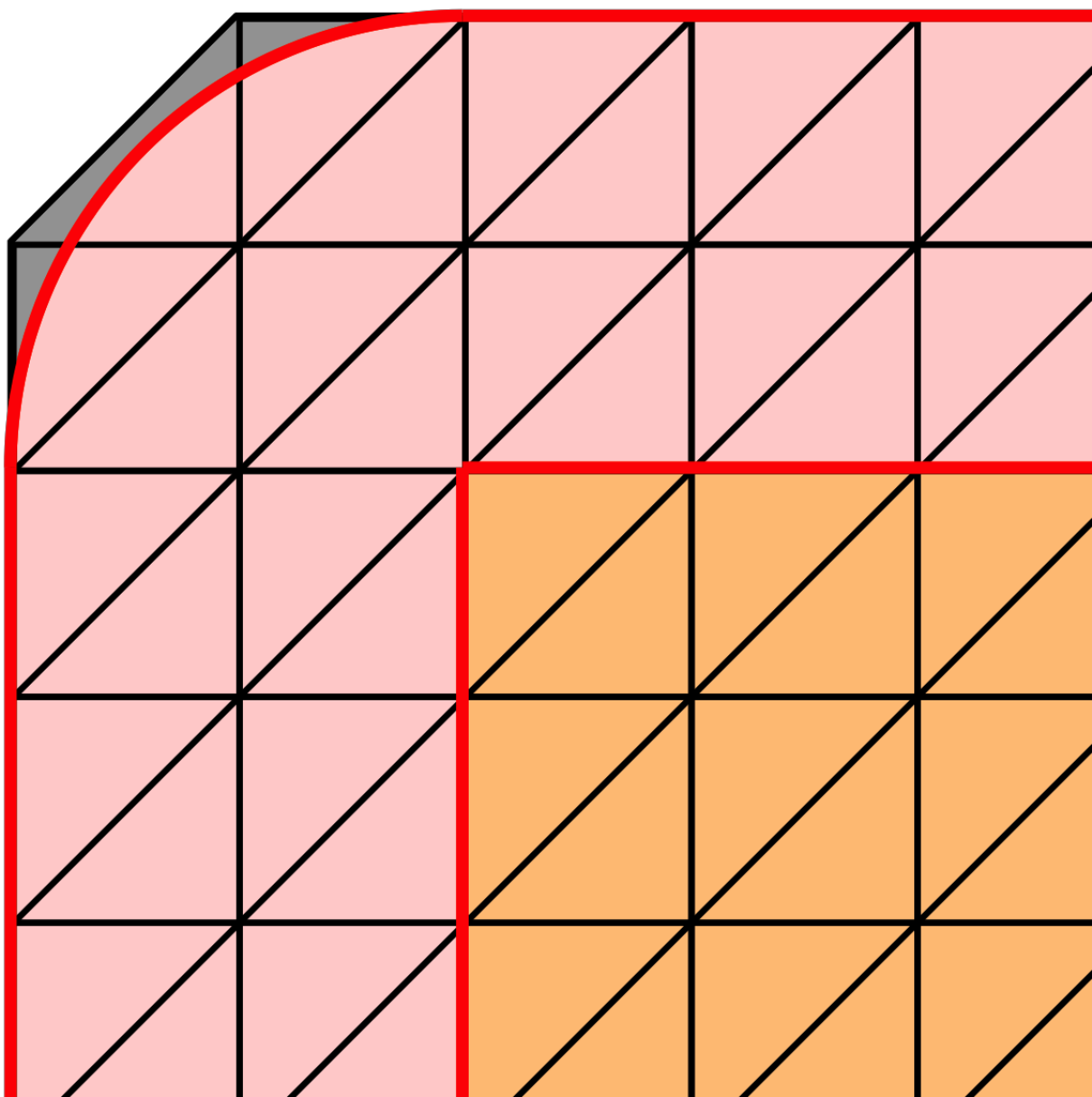}
		&
		\includegraphics[width=1.in]{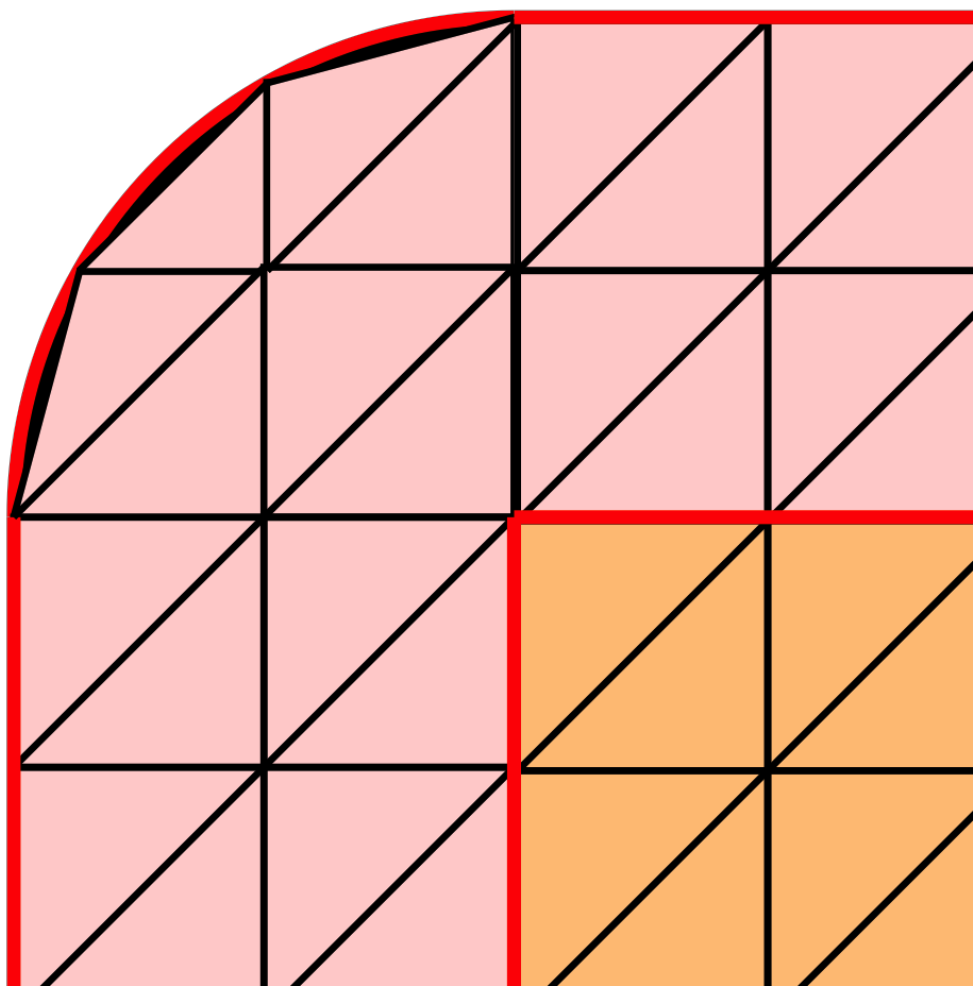}
		&
		\includegraphics[width=1.05in]{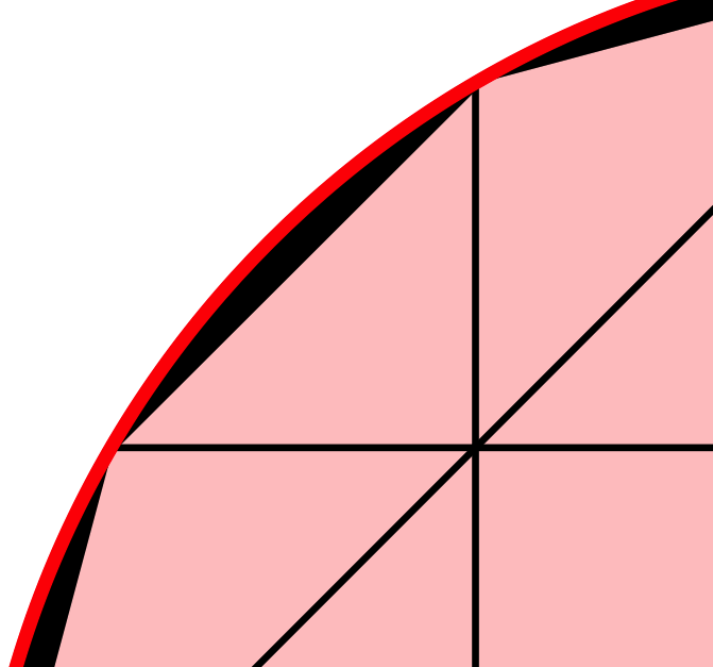}
		\\
		(a) \quad&\quad (b) \quad&\quad (c)
\end{tabular}
\caption{The pink interaction region has a curved boundary corresponding to each vertex of the orange polygonal region $\domain$ as is illustrated in (a). The curved boundary can be approximated in several ways, including using whole enclosing finite element triangles such as those that that colored in pink and gray in (a) or can be approximated by the edges of straight-sided finite elements as are illustrated in (b) and the zoom-in in (c).
}\label{fig2a}
\end{figure}

For the construction of finite element grids in the presence of curved boundaries, there exist several approaches available such as, among others, isoparametric \cite{isopar} and isogeometric \cite{isogeo} approximations. Perhaps the simplest approach is to include all triangles that intersect with the $\ell^p$-balls centered at the vertices of $\domain$. For \(p=2\) this strategy would include triangles such as those colored in both pink and gray in Figure \ref{fig2a}a. The geometric error incurred by this strategy is proportional to total area within triangles that are exterior to $\dirichletDomain$ such as the gray areas in Figure \ref{fig2a}a. That area is of ${\mathcal O}(h^2)$ for each vertex so that because, in general, there are a small number of vertices, the total geometric error induced remains of ${\mathcal O}(h^2)$; see \cite{DEliaFEM2020}. Another approach is to approximate the curved boundary by a piecewise linear approximation with each line segment defining a triangle; see Figure \ref{fig2a}b for an illustration. The geometric error incurred by this strategy is proportional to the total areas within the $\ell^2$-balls that are excluded as is illustrated by the black areas in Figures \ref{fig2a}b and \ref{fig2a}c. For each vertex, that area is of ${\mathcal O}(h^3)$ so that because, in general, there are a small number of vertices, the total geometric error induced remains of ${\mathcal O}(h^3)$; again see \cite{DEliaFEM2020}.

\begin{figure}[h!]
\centerline{
\includegraphics[width=2.5in]{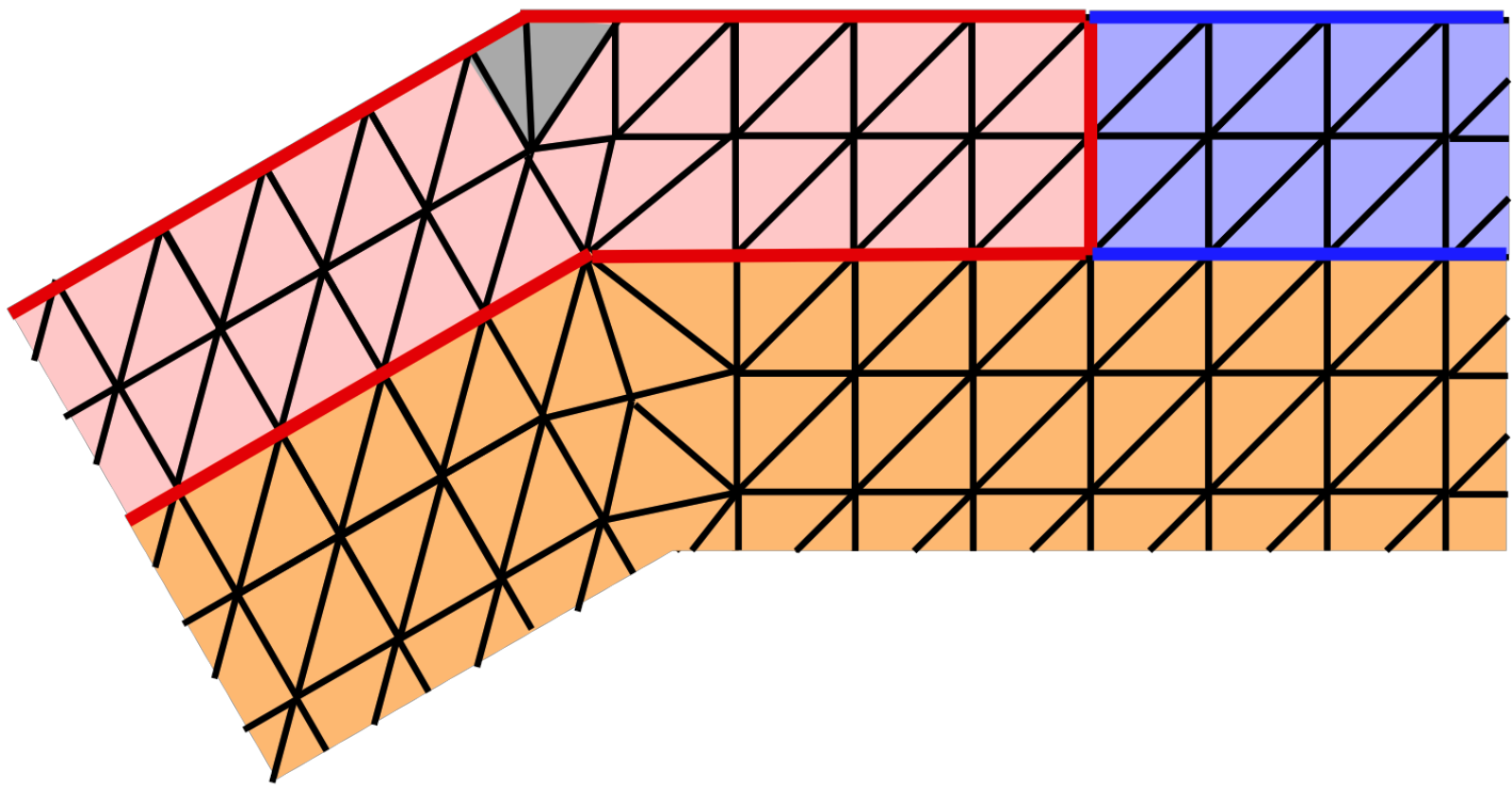}
\qquad\includegraphics[width=2.5in]{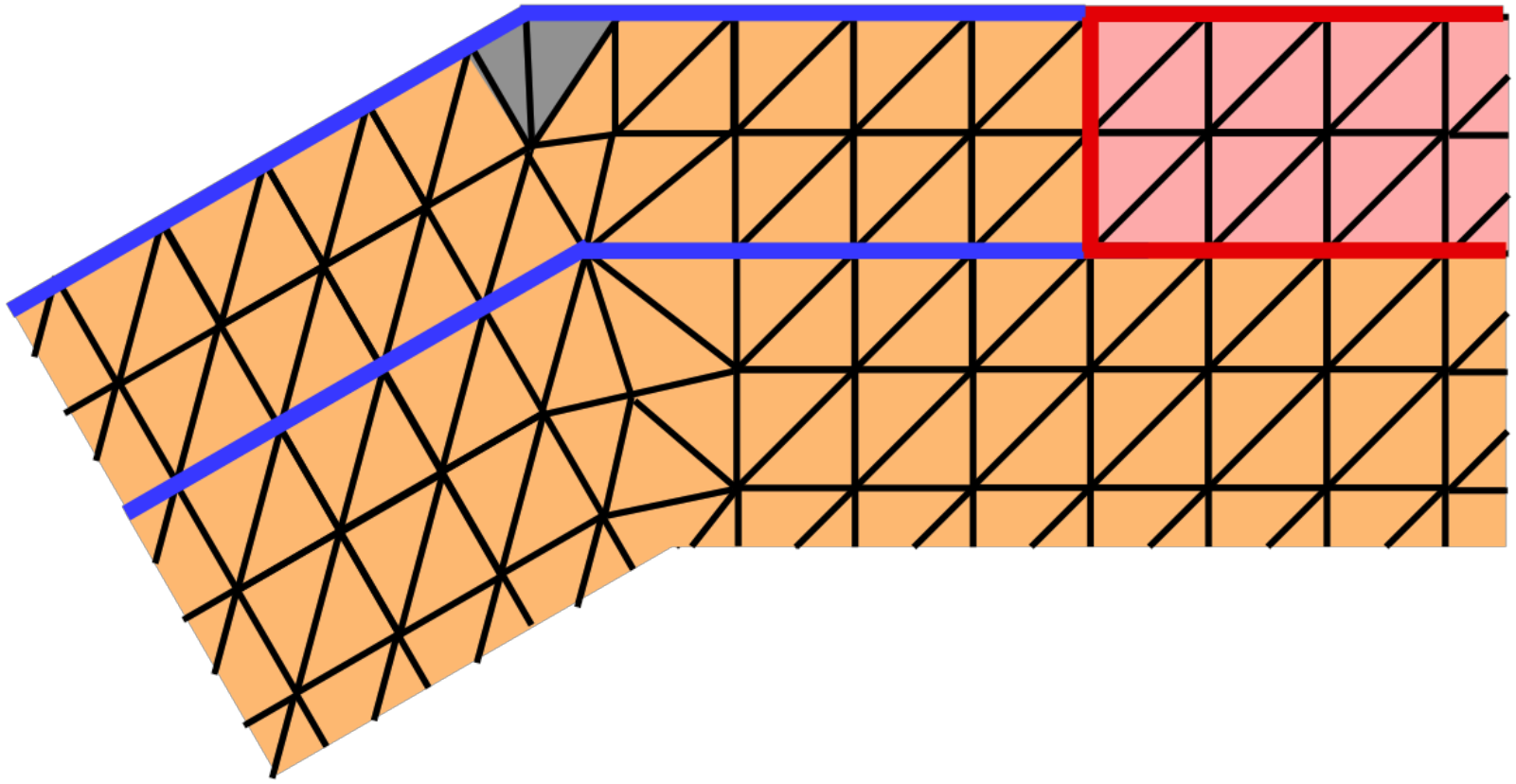}
}
\caption{
Left: A zoom in of a portion of a triangular finite element grid that respects the common boundaries of the three subregions $\domain$ (in orange), $\dirichletDomain_h$ (in pink + gray), and $\neumannDomain_h$ (in blue). Right: A triangular finite element grid with $\Omega_h=\domain\cup\Gamma^N_h$ (in orange + gray) and $\dirichletDomain_h$ (in pink) where $\Gamma_h^N$ is bounded by the blue lines. For both cases, the two gray triangles define the geometric approximation of curved boundaries induced by the use of $\ell^2$-balls in the construction of interaction regions; specifically, those triangles are of the type depicted in gray in both plots.}
\label{fig444}
\end{figure}

We denote by $\dirichletDomain_h$ and $\neumannDomain_h$ \label{inline:GammaDh} the geometric approximations of $\dirichletDomain$ and $\neumannDomain$, respectively, in case the domains have curved boundaries.
We also denote by
\begin{equation}\label{gridsgrids}
\begin{aligned}
&
\mbox{$\mcT^{\domain}$, $\mcT^D$, and $\mcT^N$ the triangulation of ${\domain}$, $\dirichletDomain_h$, and $\neumannDomain_h$, respectively}
\\[-1ex]
&
\mbox{$\mcT=\mcT^{\domain}\cup\mcT^D\cup\mcT^N$ the triangulation of $\domain\cup\dirichletDomain_h\cup\neumannDomain_h$}
\\[-1ex]
&
\mbox{$\mcT^{\Omega}=\mcT^{\domain}\cup\mcT^N$ the triangulation of $\Omega_h= \domain\cup\neumannDomain_h$.}
\end{aligned}
\end{equation}
Note that if $\neumannDomain_h=\emptyset$ we have that $\Omega_h=\Omega=\domain$. Figure \ref{fig444} illustrates the properties of the finite element grids we use, including respecting the common boundaries between subdomains and approximating vertex-induced curved boundary segments by straight-sided line segments.

\subsubsection{Approximate interaction balls}

\begin{figure}[h!]
\begin{center}
\begin{tabular}{ccc}
\includegraphics[height=1.1in]{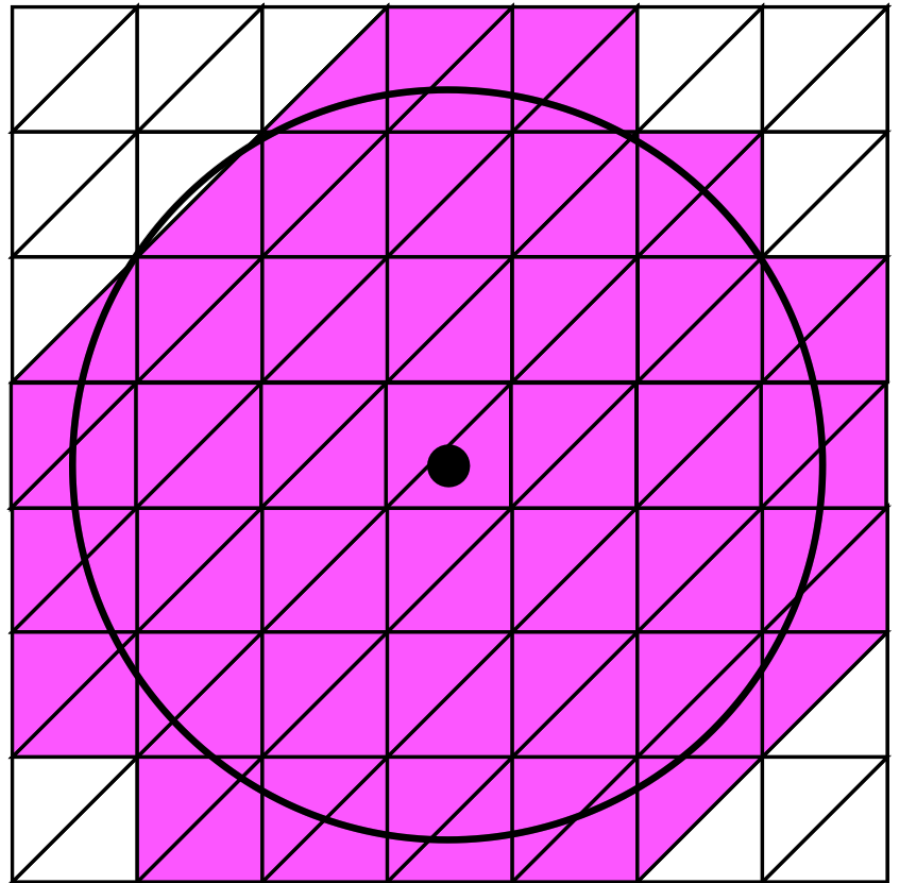}&
\includegraphics[height=1.1in]{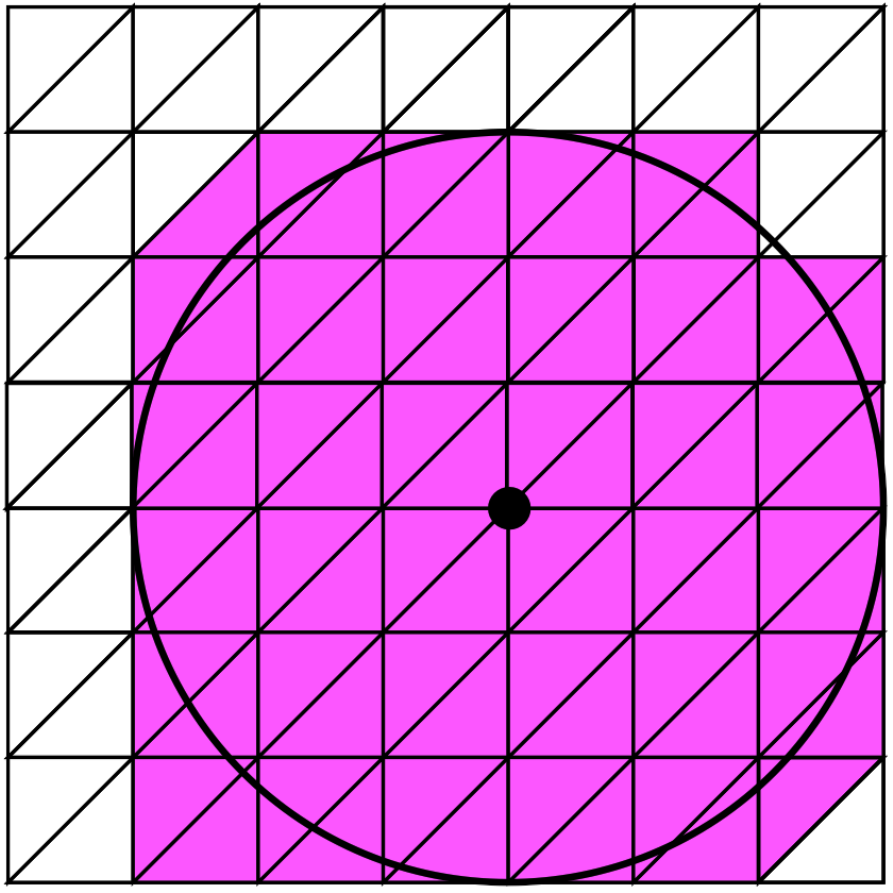}&
\includegraphics[height=1.1in]{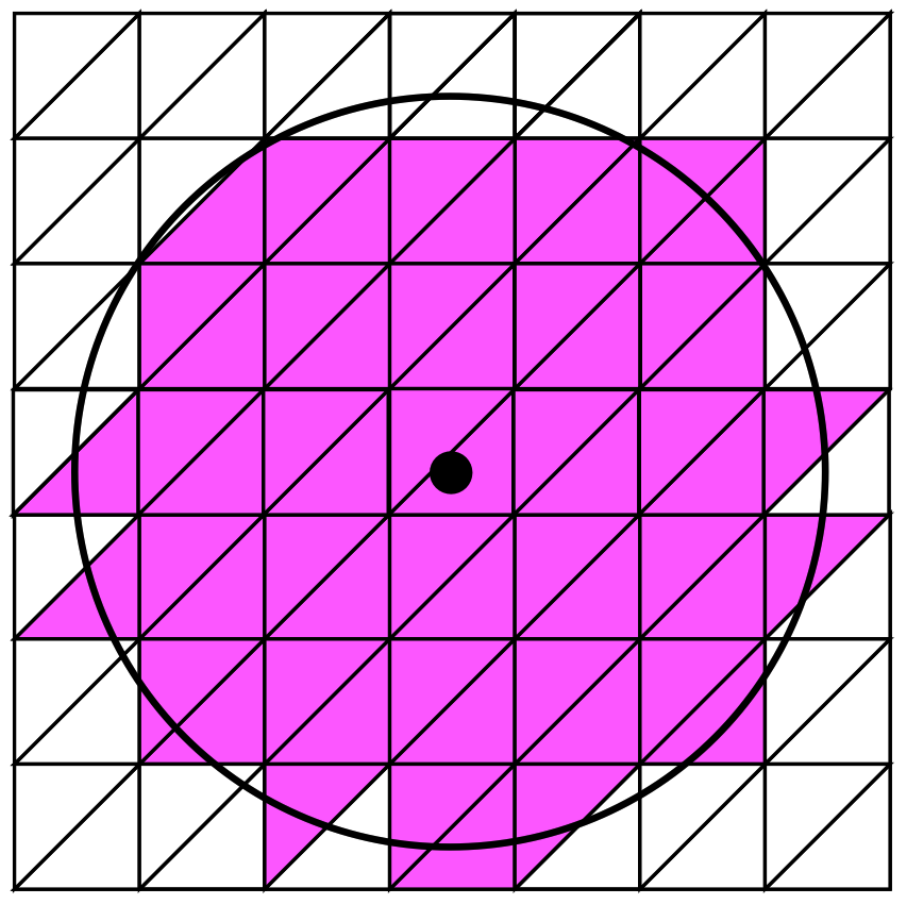}
\\
(a) $B_{\delta,h}^{123vertices}(\xb)$ &
(b) $B_{\delta,h}^{123vertices}(\xb)$ &
(c) $B_{\delta,h}^{barycenter}(\xb)$
\\[1ex]
\includegraphics[height=1.1in]{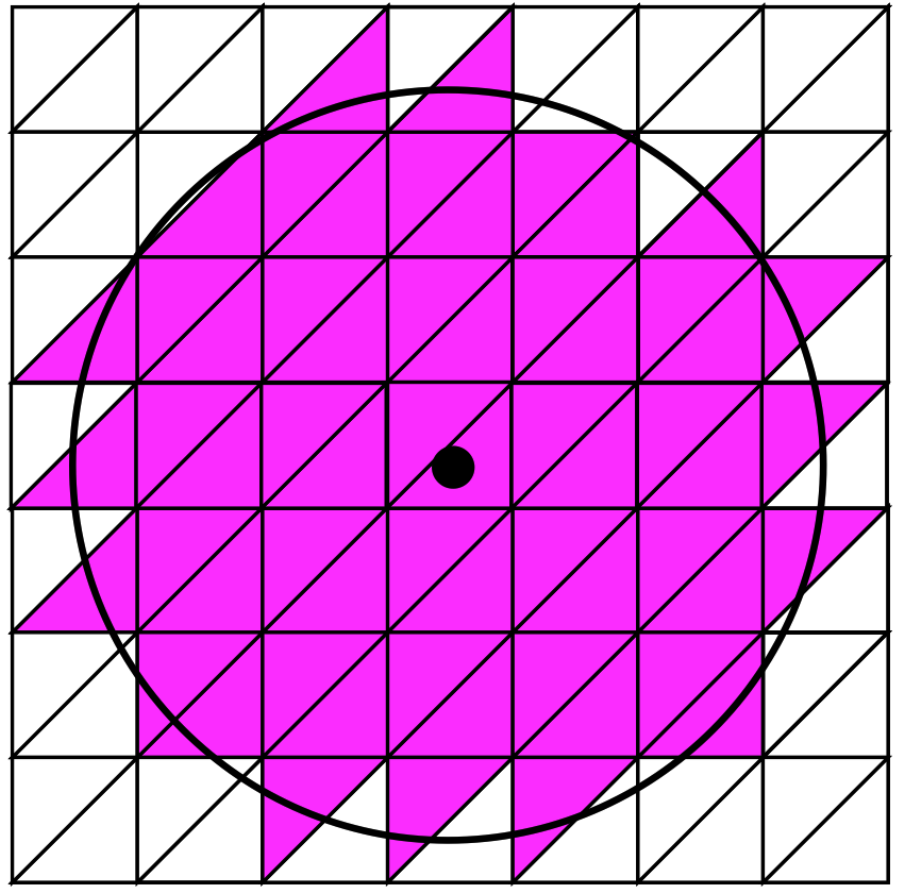}&
\includegraphics[height=1.1in]{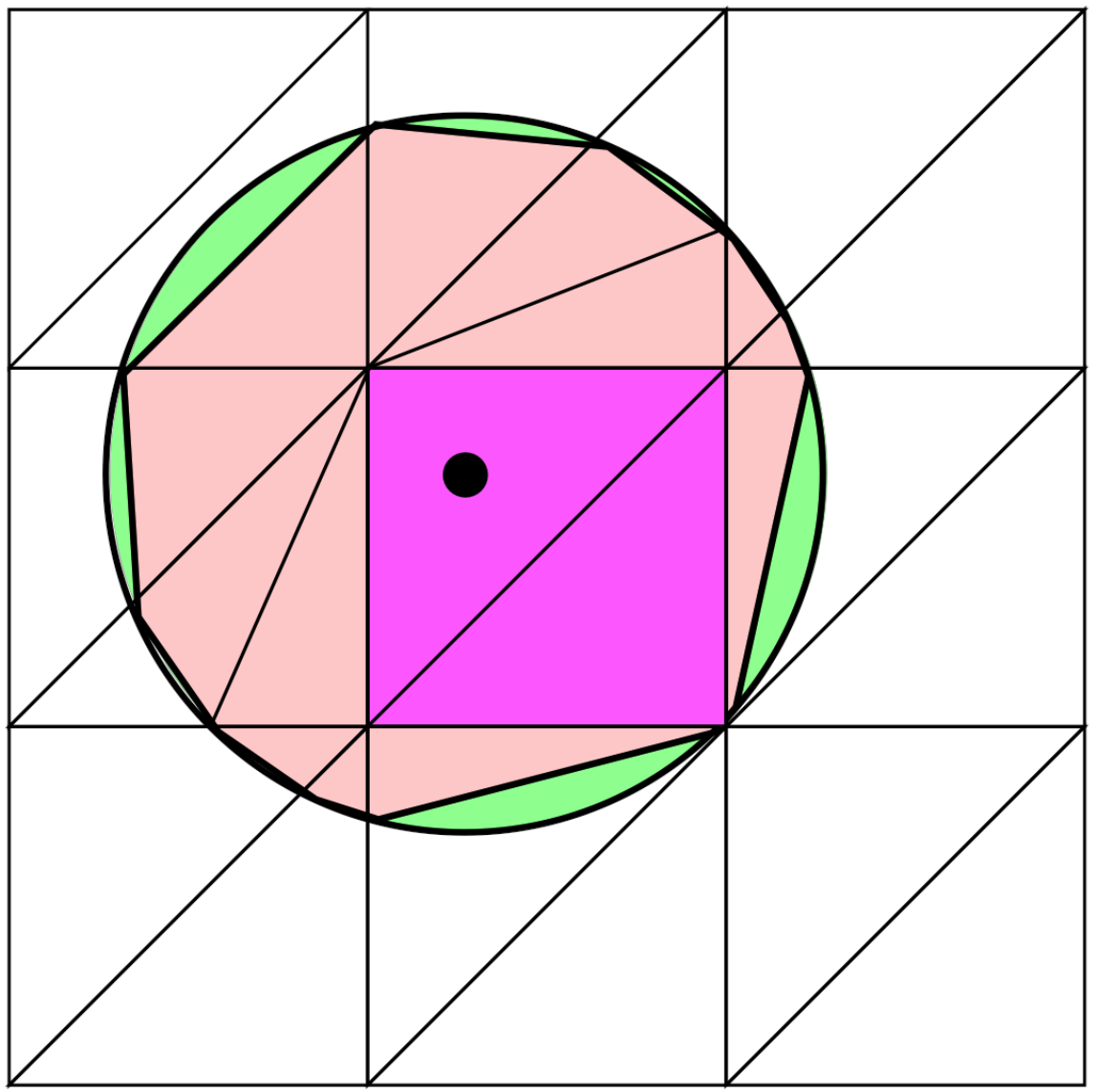}&
\includegraphics[height=1.1in]{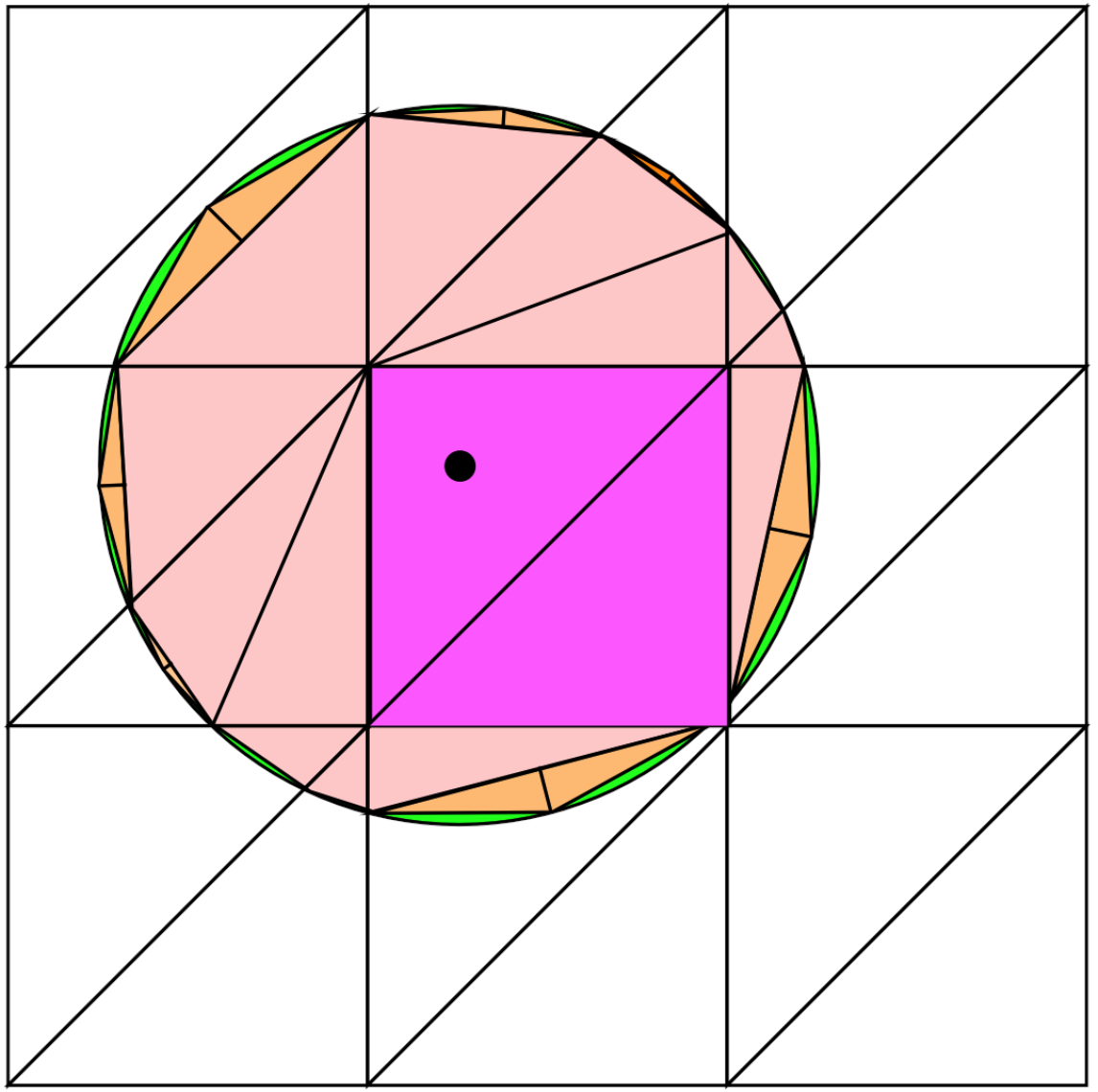}
\\
(d) $B_{\delta,h}^{23vertices}(\xb)$ &
(e) $B_{\delta,h}^{nocaps}(\xb)$ \hspace{-.17in} &
(f) $B_{\delta,h}^{approxcaps}(\xb)$ \hspace{-.17in}
\end{tabular}
\end{center}
\caption{Six approximations of the exact ball $B_{\delta,2}(\xb)$.}\label{approximate-ball}
\end{figure}

The finite element assembly of the bilinear form \(\mathcal{A}\) requires the evaluation of integrals over pairs of elements of the triangulation using quadrature.
However, special care needs to be taken to not integrate across the discontinuity of the kernel function as this can lead to slow convergence of the quadrature.
Therefore, one does not use, say, the $\ell^2$-ball to define the support of the kernel $\gamma(\xb,\yb)$. Instead most often one uses an approximate $\ell^2$-ball $B^h_{\delta,2}(\xb)$ so that instead of using $\gamma(\xb, \yb)$ as defined in \eqref{eq:scalarkernel}, we use an approximate kernel
\begin{equation}
\label{eq:scalarkernel-h}
\gamma_h(\xb, \yb) = \phi(\xb, \yb)\ind_{B^h_{\delta,2}(\xb)}(\yb) \approx \gamma(\xb, \yb).
\end{equation}
In \cite{DEliaFEM2020}, several geometric approximations the $\ell^2$-ball are discussed. For example, the first four approximate balls in Figure \ref{approximate-ball} are defined using only whole finite element triangles. For (a) and (b) in Figure \ref{approximate-ball}, we have that the approximate ball $B_{\delta,h}^{123vertices}(\xb)$ contains all finite element triangles that have at least one vertex located inside the exact ball $B_{\delta,2}(\xb)$, where the black dot denotes the center of the exact ball. Note that the location of the black dots can cause large changes in the shape of the approximate ball. For (c), the approximate ball $B_{\delta,h}^{barycenter}(\xb)$ contains all whole finite element triangles whose barycenters are located within the exact ball. For (d), the approximate ball $B_{\delta,h}^{23vertices}(\xb)$ contains all whole finite element triangles having 2 or 3 vertices located within the exact ball. Note that in (a), (c), and (d) the position of the exact ball is the same but there are large differences resulting from using different approximate balls. For (e), we have an inscribed polygon whose sides are cords of the exact ball that lie within a whole finite element triangle. Note that within this approximate ball we have whole finite element triangles (in magenta) and sub-triangles (in pink) each of which is strictly contained within a whole finite element triangle. The green caps in (e) are omitted. Finally, for (r), the green caps appearing (e) are themselves approximated by the additional sub-triangles that again are strict subsets of whole finite element triangles. As a result, the volume omitted from the exact ball is much smaller that for (e).

It will become clear in the section~\ref{sec:muultipledomain} that the choice of approximate ball impacts the subdivision of the problem into overlapping subdomains.
In what follows we will denote the evaluation of the bilinear form given in \eqref{eq:bil1} with \(\gamma\) replaced by \(\gamma_{h}\) and using numerical quadrature by \(\mathcal{A}_{h}\). \label{inline:bilQuad}
Note that \(\mathcal{A}_{h}\) only needs to be defined on the space of finite element functions, but not on the continuous space \(W(\Omega_{h}\cup\Gamma^{D}_{h})\).

\subsubsection{A discretized weak formulation}\label{fem-weak}

Given the finite element grids defined in Section \ref{fem-grid-onedomain}, we construct continuous Lagrangian, i.e., a nodal, finite element spaces of order \(k\geq 1\) as follows.
\begin{align}\label{def:femSpace}
  W_{h}(\Omega_{h}\cup\Gamma^{D}_{h})&:=\left\{v\in C(\Omega_{h}\cup\Gamma^{D}_{h},\mathbb{R}) \mid v|_{K}\in \mathbb{P}_{k}(K) ~ \forall K\in \mathcal{T} \right\}.
\end{align}
Here \(\mathbb{P}_{k}(K)\) is the space of polynomials of degree \(k\) on \(K\).
We denote by \(\{\psi_{m}^{h}\}_{m}\) the nodal basis functions spanning \(W_{h}(\Omega_{h}\cup\Gamma^{D}_{h})\) and \(\{\xb_{m}^{h}\}_{m}\) the associated nodal coordinates of the degrees of freedom.
Based on the spatial location of degrees of freedom we partition the finite element basis functions into interior and essential volume condition:
\begin{align*}
  \mathcal{I}&:= \{m \mid \xb_{m}^{h}\not\in \Gamma^{D}_{h} \}, &
  \mathcal{B}&:= \{m \mid \xb_{m}^{h}\in \Gamma^{D}_{h} \}.
\end{align*}
Let
\begin{align*}
  W_{h}^{0}(\Omega_{h},\Gamma^{D}_{h})&:=\operatorname{span}\left\{\psi_{m}^{h} \mid m\in\mathcal{I}\right\}, &
  W_{h}^{t}(\Gamma^{D}_{h})&:=\operatorname{span}\left\{\psi_{m}^{h}|_{\Gamma^{D}_{h}} \mid m\in\mathcal{B}\right\}.
\end{align*}
We denote by \(g_{h}\) an interpolation or projection of the volume data \(g\) into the finite element space \(W_{h}^{t}(\Gamma^{D}_{h})\).

A discrete weak formulation of the nonlocal volume-constrained problem \eqref{eq:weak_nonloc} is then given by
\begin{equation}\label{eq:weak_nonloc-h}
\begin{aligned}
&\mbox{\em given $f\in W'(\Omega_{h}, \Gamma^D_{h})$, $g_h(\xb)$, and $\gamma(\xb,\yb)$,
find $u_h(\xb)\in W_h(\Omega_h\cup\Gamma^{D}_{h})$ such that}
\\
&  \qquad
\mathcal{A}_{h}(u_h,v_h) = \mcF(v_h) \quad \forall\, v_h \in W^{0}_{h}(\Omega_h,\Gamma^{D}_{h})
\\&\mbox{\it and such that $u_h(\xb)=g_h(\xb)$ for $\xb\in\dirichletDomain_h$}.
\end{aligned}
\end{equation}

We denote the vector of basis coefficients of \(u_{h}\in W_h(\Omega_h\cup\Gamma^{D}_{h})\) by \(\ub_{\Omega_h\cup\Gamma^{D}_{h}}^{h}\), i.e. \(u_{h}=\sum_{m\in\mathcal{I}\cup\mathcal{B}}\ub_{\Omega_h\cup\Gamma^{D}_{h},m}^{h} \psi_{m}^{h}\).
\(\ub_{\Omega_h\cup\Gamma^{D}_{h}}^{h}\) is given by the concatenation of the vectors \(\ub^{h}\) and \(\ub_{\Gamma^{D}_{h}}^{h}\) corresponding to the coefficients of basis functions in the sets \(\mathcal{I}\) and \(\mathcal{B}\) respectively.
Similarly, we denote the coefficients of \(g_{h}\) by \(\gb^{h}\).

Therefore, since \(\ub_{\Gamma^{D}_{h}}^{h}=\gb^{h}\), \eqref{eq:weak_nonloc-h} can be equivalently rewritten as
\begin{align*}
  \mathbb{A} \ub^{h} &= \fb^{h} - \mathbb{B} \gb^{h}
\end{align*}
with
\begin{align*}
  \mathbb{A}_{m,m'}&= \mathcal{A}_{h}(\psi_{m}^{h}, \psi_{m'}^{h}) , & m,m'\in\mathcal{I},\\
  \mathbb{B}_{m,m'}&= \mathcal{A}_{h}(\psi_{m}^{h}, \psi_{m'}^{h}) , & m\in\mathcal{I}, ~m'\in\mathcal{B},\\
  \fb^{h}_{m}&= \mathcal{F}(\psi_{m}^{h}), & m\in\mathcal{I}.
\end{align*}
The matrix ${\mathbb A}$ is symmetric and if $\Gamma^D_h\neq\emptyset$, it is also positive definite.
The construction of an efficient solver for this linear system of equations is the goal of the next two sections.
We also see that non-homogeneous volume conditions enter in the form of a correction to the right-hand side vector.
For the purposes of constructing a solver we can therefore assume homogeneous volume conditions on \(\Gamma^{D}_{h}\) in what follows.

\section{Nonlocal models for diffusion: the multi-domain domain decomposition setting}\label{sec:muultipledomain}

\subsection{Nonlocal subdivision}

In the partial differential equation setting, substructuring-based domain decomposition methods partition the domain into covering, non-overlapping subdomains. This is not viable for nonlocal models because nonlocality necessarily induces overlaps between subdomains, i.e., the interaction region for a subdomain is in part constituted of finite element triangles belonging to abutting subdomains.

\begin{remark}
{\em
In this section, we define several geometric entities that are then used in Section \ref{sec:mdwf} to define the subdomain problems. How those geometric entities are constructed in practice is discussed in Section \ref{sec:implementation}.
}\qquad$\Box$
\end{remark}

For the building of substructures (which we also refer to as subdomains), we are given single-domain finite element grids $\mcT^{\Omega_h}$ and $\mcT^{\dirichletDomain_h}$ covering $\Omega_h$ and $\dirichletDomain_h$, respectively, into finite elements, e.g., as discussed in Section \ref{fem-grid-onedomain} and defined in \ref{gridsgrids}.

\begin{figure}
	\centering
	\begin{subfigure}{0.49\linewidth}
		\includegraphics[page=1,width=1.\linewidth]{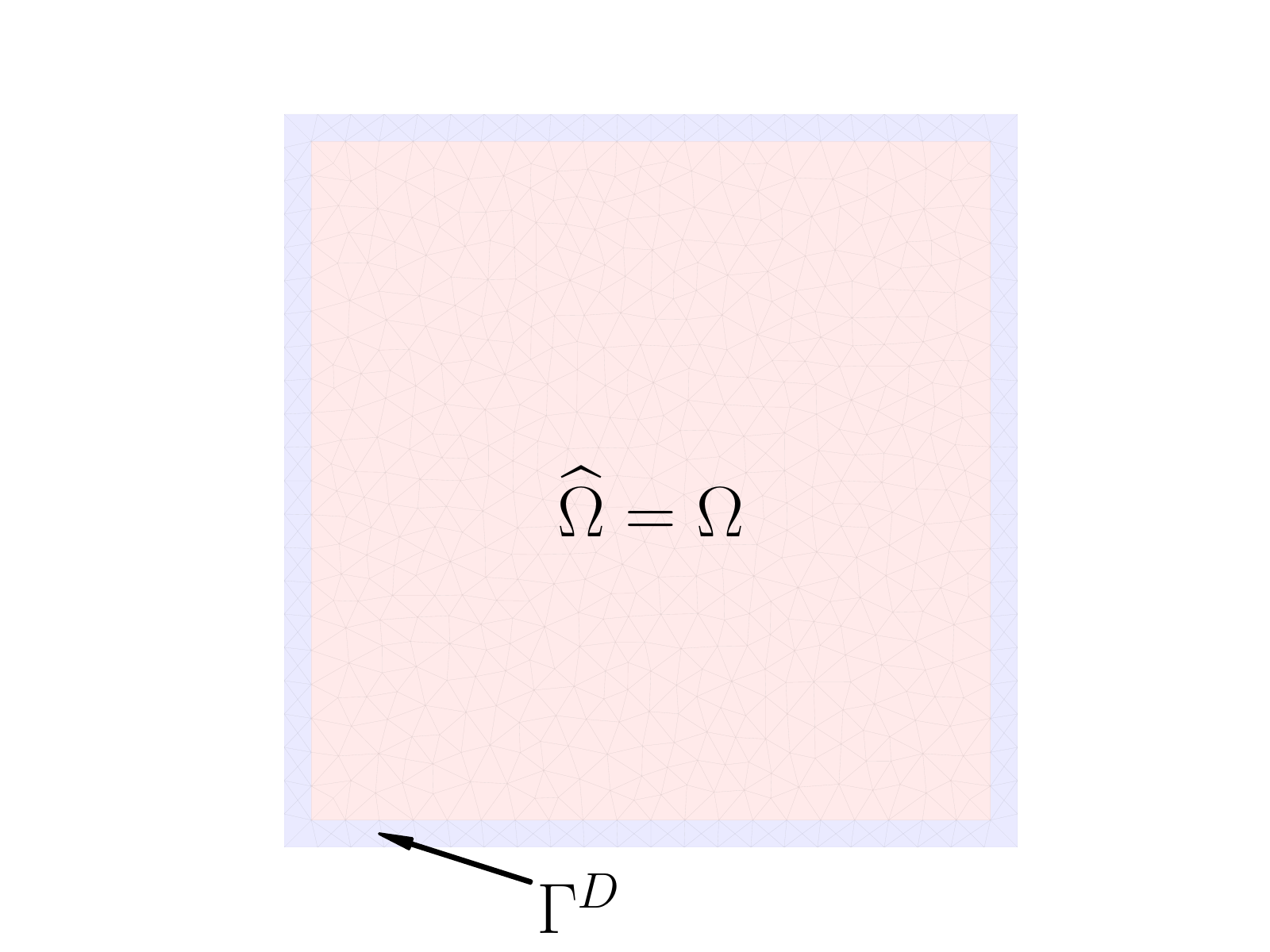}
		\caption{}
	\end{subfigure}
	\begin{subfigure}{0.49\linewidth}
		\includegraphics[page=2,width=1.\linewidth]{Figures/unregular_subdivision}
		\caption{}
	\end{subfigure}
	\caption{
		(a) The single domain $\widehat{\Omega}$ (red) and the Dirichlet domain $\Gamma^D$ (blue).
		(b) Decomposition of ${\Omega}$ into 9 substructures.
		The darkness of the colors indicates overlaps.
		The domains $\Omega_1, \dots, \Omega_9$, defined in \eqref{eq:subdomainDefinitions},
		do not touch.
	}
	\label{fig:mesh-decomposition}
\end{figure}
\begin{figure}
	\centering
	\begin{subfigure}{0.49\linewidth}
		\includegraphics[page=3,width=1.\linewidth]{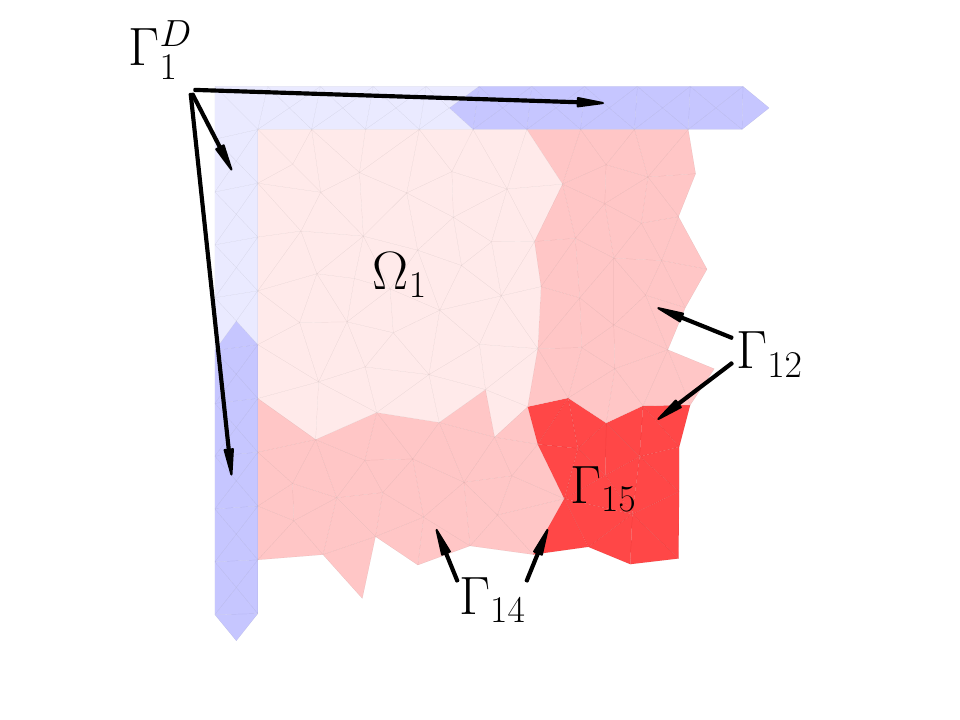}
		\caption{}
	\end{subfigure}
	\begin{subfigure}{0.49\linewidth}
		\includegraphics[page=4,width=1.\linewidth]{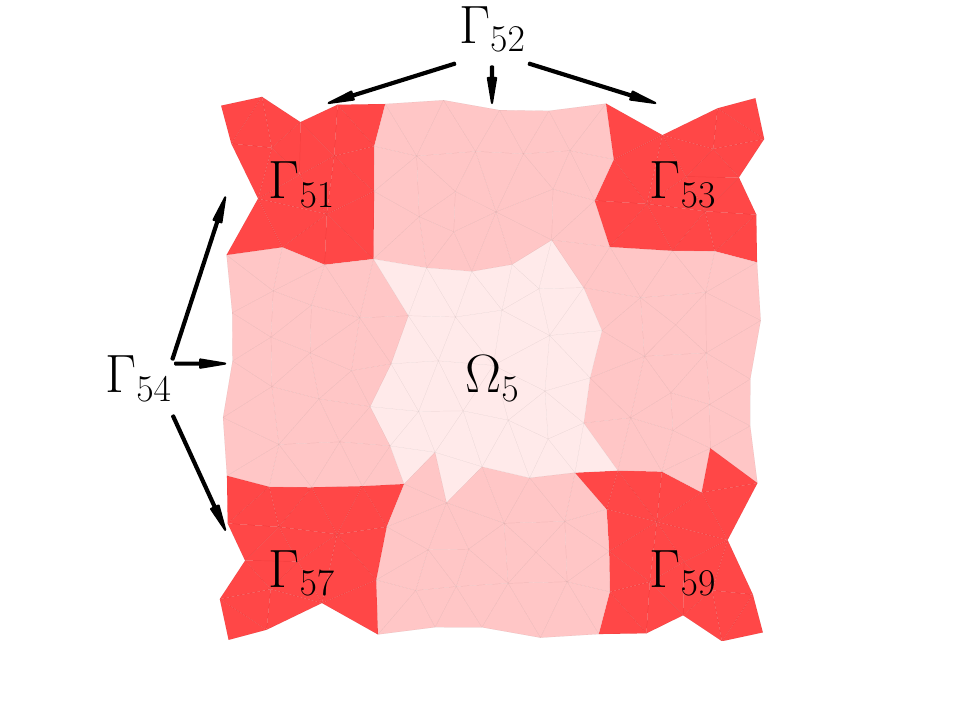}
		\caption{}
	\end{subfigure}
	\caption{
		(a) The overlap of $\widehat\Omega_1 \cup \Gamma^D_1$ with it's neighboring subdomains contains
		$\Gamma_{12}$, $\Gamma_{14}$, and $\Gamma_{15}$. The corresponding Dirichlet domain
		is given by $\Gamma_1^D$ \eqref{eq:subdomainDefinitions}.
		(b) The floating subdomain
		$\widehat\Omega_5 = \Omega_5 \cup \Gamma_5$ has an empty intersection with $\Gamma^D$.
		The dark red parts overlap with 3 other subdomains. See for example $\Gamma_{51}$
		which has overlaps with
		$\widehat\Omega_1$,
		$\widehat\Omega_2$ and
		$\widehat\Omega_4$.
	}
	\label{fig:subdomains}
\end{figure}

While the construction of a provably viable subdivision of $\Omega_{h} \cup \Gamma^D_{h}$ can be an intricate task, a nonconstructive definition is easily obtained.
We introduce a subdivision of $\Omega_{h} \cup \Gamma^D_{h}$ which is
necessary and sufficient for an equivalence of the single and (yet to be introduced) multidomain formulation of \eqref{eq:weak_nonloc}.
To that end, we define the support of a kernel by
\begin{align}
\label{eq:suppKernel}
\support(\gamma_{h}) = \overline{\lbrace (\xb, \yb)\in (\Omega_{h} \cup \Gamma^D_{h})^{2} ~\colon ~ \gamma_{h}(\xb, \yb)  \neq 0 \rbrace}
\approx \overline{\lbrace (\xb, \yb)\in (\Omega_{h} \cup \Gamma^D_{h})^{2} ~\colon ~ |\xb - \yb|_p \leq \delta\rbrace}.
\end{align}
\begin{definition}
	\label{def:nonlocalsubdivision}
	A \textit{nonlocal subdivision} $\family{\widehat{\Omega}_{k}\cup\Gamma^D_{k}}$ of $\Omega_{h} \cup \Gamma^D_{h}$ is a family of subsets $\widehat{\Omega}_{k}\cup\Gamma^D_{k}$ with
	$\widehat{\Omega}_k \subset \Omega_{h}$ and
	$\Gamma^D_k \subset \Gamma^D_{h}$ for $k=1,\dots,K$
	that fulfill
	\begin{align}
	\label{eq:subRequirement}
	\support(\gamma_{h})
	\cap 
	\left[(\Omega_{h} \cup \Gamma^D_{h}) \times (\Omega_{h} \cup \Gamma^D_{h})\right]
	\subset 
	\bigcup_{k=1}^{K} 
	(\widehat{\Omega}_{k}\cup\Gamma^D_{k})
	\times 
	(\widehat{\Omega}_{k}\cup\Gamma^D_{k}).
	\end{align} 
\end{definition}

Definition \ref{def:nonlocalsubdivision} guarantees that the collection of subproblems defined on $\family{\widehat{\Omega}_{k}\cup\Gamma^D_{k}}$ can recover all nonlocal interactions which appear in the single domain. 
The substructures $\widehat{\Omega}_{k}\cup\Gamma^D_{k}$ can chosen to be larger than necessary.
This leads to redundant computations and more communication between subdomain than strictly necessary but it also allows to use coverings which do not cut through finite elements, see Figure \ref{fig:mesh-decomposition}.
Based on the substructures $\widehat{\Omega}_{k}\cup\Gamma^D_{k}$ we define 
\begin{equation}
\label{eq:subdomainDefinitions}
\begin{aligned}
\Gamma_{k \ell} &:= \widehat{\Omega}_k \cap \widehat{\Omega}_\ell \textrm{ for } k \neq \ell,	\\
\Gamma_k &:= \bigcup_{\ell=1, \ell\neq k}^K \Gamma_{k \ell}, \\
\Omega_k &:= \widehat{\Omega}_k \setminus \Gamma_k, \\
\Gamma &:= \bigcup_{k=1}^K \Gamma_k.
\end{aligned}
\end{equation} 
If the substructures $\widehat{\Omega}_{k}\cup\Gamma^D_{k}$ do not cut through finite elements, then  
the sets $\Gamma_{k \ell}, \Gamma_k$ and $\Omega_k$ do not cut through elements either.
Note
that $\Gamma^D_k$ and $\Gamma_{k \ell}$ are possibly empty. 
If $\Gamma^D_k$ is empty $\widehat{\Omega}_k = \Omega_k \cup \Gamma_k$ is called a \emph{floating subdomain}, see Figure \ref{fig:subdomains}.
However, $\Gamma_k \cup \Gamma^D_k$ has positive measure for every \(k\).
The set $\Gamma_k$ is the part of the interface  $\Gamma$ interacting with the substructure $\widehat{\Omega}_{k}\cup\Gamma^D_{k}$.
It is not a subset of the region $\Gamma^N$ where Neumann boundary conditions are imposed on the single domain problem,
but it is indeed the region where Neumann conditions are imposed on the subproblem. 
Ultimately, the multi-domain system can be reduced to a problem on $\family{\Gamma_k}$,
as we will see in Section \ref{subsec:DiscretizationOfSubdomains}. 
\begin{remark}
The overlaps of the domains $\Gamma_{k\ell}$ in Figure \ref{subsec:DiscretizationOfSubdomains} are smaller than the overlap of the Dirichlet boundaries $\Gamma_{k}^D \cap \Gamma_{\ell}^D$. This is a result to the specific construction algorithm we apply to obtain the decomposition. We shortly outline it in Section \ref{subsec:buildingSubtructures} and refer the reader to \cite{capodaglio2020general} for a thorough construction of a nonlocal decomposition.
\end{remark}
\subsection{Multi-domain weak formulations}\label{sec:mdwf}

Again we restrict the presentation to the nonlocal diffusion model based on the kernel $\gamma_h(\xb,\yb)=\phi_h(\xb,\yb)\ind_{B_{\delta,p}^{h}(\xb)}(\yb)$; the nonlocal solid mechanics case is treated in an entirely analogous manner.

Nonlocal subdivisions are necessarily overlapping yet they must avoid multivalued solutions in regions of overlaps. Thus, we have to impose the  constraint
\begin{equation}\label{Gamma_kGamma_k'}
     u_k(\xb) = u_{\ell}(\xb)\quad\mbox{whenever \,\,\, $\Gamma_{k\ell} \neq \emptyset$\quad for $k,\ell=1,\ldots,K$}.
\end{equation}
The overlap also causes the necessity of introducing weights that ensure that the duplications of terms in the bilinear forms do not occur. To this end, following \cite{capodaglio2020general}, we define the function
\begin{equation}
\label{eq:zeta}
\zeta(\xb, \yb) = 	\sum_{\kSub=1}^\KSub
\mcX_{\widehat{\Omega}_{k}\cup\Gamma^D_{k}}(\xb)
\mcX_{\widehat{\Omega}_{k}\cup\Gamma^D_{k}}(\yb)
\end{equation}
which counts the number of overlaps of the subdomains.
The integer-valued counting function $\zeta(\xb, \yb)$ is symmetric, piecewise constant by element in \(\mathcal{T}\) and $0 \leq \zeta(\xb, \yb) \leq K$ for all $\xb, \yb \in \Omega \cup \Gamma^D$. 
In general $\zeta(\xb,\yb)$ can attain $0$ for some pairs of (remote) points, however due to \eqref{eq:subRequirement} and Definition \ref{def:nonlocalsubdivision}  we have that $1 \leq \zeta(\xb, \yb)$ on the support of the kernel $\gamma(\xb,\yb)$ and $\family{\zeta^{-1}(\xb,\yb) \mcX_{\widehat{\Omega}_{k}\cup\Gamma^D_{k}}(\xb) \mcX_{\widehat{\Omega}_{k}\cup\Gamma^D_{k}}(\yb)}$ defines a partition of unity on that support. Hence, the weighted kernel $\zeta^{-1}(\xb,\yb) \gamma(\xb,\yb)$ is well defined.

Corresponding to each of the $K$ subdomains $\widehat{\Omega}_k\cup\Gamma_k$, we define, for $k=1,\ldots,K$, the bilinear forms
\begin{align}\label{eq:bilinear-multi1}
  {\mathcal A}_{k}(u_k,v_k) &:=
                            \int_{\widehat{\Omega}_{k}\cup\Gamma^{D}_{k}} \int_{\widehat{\Omega}_{k}\cup\Gamma^{D}_{k}} \big(u_k(\yb) - u_k(\xb)\big) \big(v_k(\yb) - v_k(\xb)\big) \gamma(\xb,\yb)\zeta(\xb,\yb)^{-1}  d\yb  d\xb.
\end{align}
and let \(\mathcal{A}_{k,h}\) be the evaluation of \(\mathcal{A}_{k}\) using numerical quadrature and the kernel approximation \(\gamma_{h}\).
Similarly, we can define the linear forms
\begin{align}\label{eq:bilinear-multi5}
  \mcF_{k}(v_k) &= \int_{{\widehat{\Omega}_{k}}} f(\xb) v_{k}(\xb) \zeta(\xb,\xb)^{-1} d\xb.
\end{align}

\begin{remark}\label{rem:floating}
  {\em
    For floating subdomain $\Gamma^D_k$ is empty and no Dirichlet condition is enforced.
    As a result ${\mathcal A}_k(u_k,v_k)=0$ if $u_k(\xb)=$ {\em constant}, implying that the subdomain problem is singular.
  }
\qquad$\Box$
\end{remark}

We can then define $K$ subdomain weak formulations $\mcA_k(u_k,v_k)=\mcF_{k}(v_k)$ for $k=1,\ldots,K$.
Of course, this needs to be supplemented by the constraints in \eqref{Gamma_kGamma_k'} and volume conditions.
Using the subdomain energy functionals
\begin{align}
\label{eq:subE}
\mcE_k(u_k) = \frac{1}{2} \mcA_{k,h}(u_k, u_k) - \mcF_k(u_k)
\end{align}
allows to state the multi-domain formulation \cite{capodaglio2020general} of the nonlocal diffusion problem as
\begin{equation}\label{eq:subWeakNonlocal}
\begin{aligned} 
&\mbox{given $f\in
	W'(\Omega_{h}, \Gamma^D_{h})$, }
\\
&\mbox{and a nonlocal decomposition $\family{\widehat{\Omega}_{k}\cup\Gamma^D_{k}}$}, \\
&\mbox{find $u_k \in W_{h}(\widehat{\Omega}_{k}\cup\Gamma^D_{k})$ such that}\\
&  \qquad
\sum_{k=1}^{K}\mathcal{E}_k(u_k) = 
\inf\limits_{
	\substack{
		v_k\in W(\widehat{\Omega}_{k}\cup\Gamma^D_{k}),\\ 
		k=1,\dots,K
	}
} 
\sum_{k=1}^{K}\mathcal{E}_k(v_k) \\
&
\mbox{subject to, $u_k(\xb)=0$ for $\xb\in\Gamma^D_k$,}\\
&
\mbox{~~ $u_k(\xb)=u_j(\xb)$ for $k \neq \ell$ and all $\xb\in \Gamma_{k \ell}$}.
\end{aligned}
\end{equation}
Here of course $\Gamma_{k \ell}$ is empty if 
$\widehat{\Omega}_k$ and 
$\widehat{\Omega}_\ell$ are disjoint.
The equivalence of \eqref{eq:subWeakNonlocal} and \eqref{eq:weak_nonloc} is shown in \cite[Proposition 9]{capodaglio2020general}.

\subsection{Matrix formulation of the multi-domain problem}
\label{subsec:DiscretizationOfSubdomains}
Consider again Figure \ref{fig:mesh-decomposition}. The domains $\Omega_k$ interact  via \(\gamma_{h}\) only with points inside $\widehat{\Omega}_{k}\cup\Gamma^D_{k}$.
Hence, after a solution has been determined on $\Gamma$,
the unknowns associated  to $\Omega_k$ can be determined in parallel, without communication among the subdomains.
In this section we therefore split the discretization of \eqref{eq:subWeakNonlocal} into a constrained problem with unknowns associated to $\Gamma_k$  and a backward substitution which determines the solutions on $\Omega_k$.
This requires to split the energy functional $\mcE_k(\cdot)$ by a block decomposition which involves a Schur complement.

Similar to the splitting of the finite element space in section~\ref{fem-weak} we introduce a splitting of \(W^{0}_{h}(\widehat{\Omega}_{k},\Gamma^{D}_{k})\) into unknowns that do not interact with other subdomains and unknowns that do interact with other subdomains.
\begin{align*}
  \mathcal{I}_{k}&:= \{m \mid \xb_{m}^{h}\in \Omega_{k} \}, &
  \mathcal{G}_{k}&:= \{m \mid \xb_{m}^{h}\in \Gamma_{k} \}
\end{align*}
and let
\begin{align*}
  W_{h}^{0}(\Omega_{k},\Gamma_{k}\cup\Gamma^{D}_{k})&:=\operatorname{span}\left\{\psi_{m}^{h} \mid m\in\mathcal{I}_{k}\right\}, &
  W_{h}(\Gamma_{k})&:=\operatorname{span}\left\{\psi_{m}^{h} \mid m\in\mathcal{G}_{k}\right\}.
\end{align*}

As in the single domain setting before, we can identify a finite element function \(u_{k}\in W^{0}_{h}(\widehat{\Omega}_{k},\Gamma^{D}_{k})\) with its coefficient vectors:
\begin{align*}
  u_{k}&= \sum_{m\in\mathcal{I}_{k}}\ub_{\Omega,m}^{h,k}\psi_{m}^{h} + \sum_{m\in\mathcal{G}_{k}}\ub_{\Gamma,m}^{h,k}\psi_{m}^{h}.
\end{align*}
Here, we have already split into coefficients \(\ub_{\Omega}^{h,k}\) and \(\ub_{\Gamma}^{h,k}\) associate with \(\Omega_{k}\) and \(\Gamma_{k}\) respectively.
We denote $U_{\Gamma,k} := \R^{M_k}$, $M_k := |\mathcal{G}_k|$, the spaces of coefficients corresponding to $W_h(\Gamma_k)$,
and $U_{\Omega,k}:=\R^{|\mathcal{I}_{k}|}$ which corresponds to  $W^{0}_h(\Omega_k,\Gamma_{k}\cup\Gamma^{D}_{k})$.
The discrete reduced multi-domain problem seeks a solution \(\ub^h = \family{\ub_\Gamma^{h,k}}\) in the product space $U := \prod_{k=1}^K U_{\Gamma, k}$. \label{inline:U}

It is natural to restrict ourselves to the space $U$ when we formulate the constraints
because in \eqref{eq:subWeakNonlocal} there are no constraints imposed on points in $\Omega_k$.
Of course $U$ 
contains elements which are \textit{infeasible} in \eqref{eq:subWeakNonlocal}.
The subspace of feasible coefficients is expressed as the null-space of
a matrix  $(\Bb^1, \dots, \Bb^K) = \Bb$ with entries $\lbrace 0, 1, -1\rbrace$.
More precisely, we have for some $\ub^h \in U$ that $\Bb(\ub^h) = 0$ if and only if
$(\ub^h_k)_{i_k}=(\ub^h_\ell)_{i_\ell}$ 
for all $k \neq \ell$ with $\Gamma_{k \ell} \neq \emptyset$
and all $\xb^h_{i_k} = \xb^h_{i_\ell}$ for $\xb^h_{i_k} \in \Gamma_k$ and $\xb^h_{i_\ell} \in \Gamma_{\ell}$.
We find that for each $\xb^h \in \Gamma$ there are $\zeta(\xb^h, \xb^h)$ corresponding coefficient values in $\ub^h$ which leads to at least $\zeta(\xb^h, \xb^h)-1$ constraints per point.
\begin{remark}
\label{rem:zetajumps}
    The function $\zeta(\xb, \yb)$ counts the overlaps of the subdomains.
    While it is constant on the (open) elements itself it might exhibit 
    jumps along the element boundaries. However, it is always clear whether a vertex $\xb^h$ belongs to a subdomain $\Omega_k$ or not and the number of overlaps in that point can be assessed. Therefore, $\zeta(\xb^h, \xb^h)$ is well-defined.
\end{remark}
If we assume that the kernel is scalar-valued and the matrix $\Bb$ has full row-rank, so that there are no redundancies, we obtain
for
\begin{align}
\label{eq:constraintSpace}
\Lambda := \R^{M_C} \textrm{ with } M_{C} := \sum_{\xb^h \in \Gamma^h} (\zeta(\xb^h, \xb^h) -1)
\end{align}
that  $\range(\Bb) = \Lambda$. For tensorial kernels we require constraints for each component of the ansatz functions.
Note that the definition of $\Bb$ is not unique.

We now turn to the discretization of the subdomain bilinear forms \({\mathcal A}_{k,h}\) on \(W^{0}_{h}(\widehat{\Omega}_{k},\Gamma^{D}_{k})\).
The corresponding discrete stiffness matrices are given by
\begin{align}
  \label{eq:subStiffnessMatrix}
&(\mbA^{k})_{m \ell} := \mcA_{k,h}(\psi_m^{h}, \psi_\ell^{h})
\textrm{ for } m,\ell \in \mathcal{I}_{k}\cup\mathcal{G}_{k}.
\end{align}
Their assembly can be performed in a distributed manner and we give more details in Section \ref{sec:implementation}.
\(\mbA^{k}\) can be given in block form in terms of the four submatrices
\begin{align}
\label{eq:subStiffnessMatrixSubBlocks}
\mbA^{k} = \begin{pmatrix}
\mbA_{\Omega \Omega}^{k} 		& \mbA_{\Omega \Gamma}^{k} \\
\mbA_{\Gamma \Omega}^{k}	& \mbA_{\Gamma \Gamma}^{k}
\end{pmatrix}.
\end{align}
corresponding to indices \(\mathcal{I}_{k}\times\mathcal{I}_{k}\), \(\mathcal{I}_{k}\times\mathcal{G}_{k}\), \(\mathcal{G}_{k}\times\mathcal{I}_{k}\) and \(\mathcal{G}_{k}\times\mathcal{G}_{k}\).

While a reduction of the constraint matrix to unknowns in $U$ is completely natural, 
we require an application of a block decomposition to achieve the same for the objective in \eqref{eq:subWeakNonlocal}.

The space $W^{0}_h(\Omega_k,\Gamma_k \cup \Gamma^D_k) \subset W^0(\Omega_k, \Gamma_k \cup \Gamma^D_k)$ and $\mbA_{\Omega \Omega}^{k}$ is non-singular
as $\Gamma_k \cup \Gamma^D_k$ has nonzero measure. 
This is true in particular  
if $\Gamma^D_k$ is empty and
$\widehat{\Omega}_k$
is a floating subdomain in which case the matrix $\mbA^k$ itself is not invertible.
As $\mbA^k$  is symmetric and $\mbA_{\Omega \Omega}^{k}$ is regular we can represent it by the block Cholesky decomposition 
\begin{align}
\label{eq:splittingBlockCholesky}
\mbA^k = \begin{pmatrix}
I							 								& 0 \\
\mbA_{\Gamma \Omega}^{k} (\mbA _{\Omega \Omega}^k)^{-1}	& I
\end{pmatrix}
\begin{pmatrix}
\mbA_{\Omega \Omega}^{k} 		& 0 \\
0									& \mbS_k
\end{pmatrix}
\begin{pmatrix}
I									& (\mbA _{\Omega \Omega}^k)^{-1} \mbA_{\Omega \Gamma}^{k} \\
0									& I
\end{pmatrix},
\end{align}
where
\begin{align}
\label{eq:subSchurComplement}
\mbS^{k} =  \mbA_{\Gamma \Gamma}^{k} - \mbA_{\Gamma \Omega}^{k} (\mbA_{\Omega \Omega}^{k})^{-1} \mbA_{\Omega \Gamma}^{k}
\end{align}
is the Schur complement.
Note that there is no need to store $\mbS^k$ as it is only required in terms of matrix vector products.
The discretization $\tilde{\fb}^{h,k}_{\Gamma} \in U_{\Gamma,k}$ and 
$\tilde{\fb}^{h,k}_{\Omega} \in U_{\Omega,k}$
of the linear functional \eqref{eq:bilinear-multi5} is given by
\begin{align*}
& (\tilde{\fb}^{h,k}_{\Omega})_m = \mcF_k(\psi_m^h) \textrm{ for } m\in\mathcal{I}_{k}, \\
& (\tilde{\fb}^{h,k}_{\Gamma})_m = \mcF_k(\psi_m^h) \textrm{ for } m\in\mathcal{G}_{k}.
\end{align*}
Now, let \(u_{k}\) be given in terms of the coefficient vectors $(\ub_\Omega^{h,k}, \ub_\Gamma^{h,k}) \in U_{\Omega,k} \times U_{\Gamma,k}$
and define $\vb^{h,k}_\Omega \in U_{\Omega, k}$ by
\begin{align}
\label{eq:splittingBasisChange}
\begin{pmatrix}
\vb^{h,k}_\Omega \\
\ub^{h,k}_\Gamma
\end{pmatrix}
=
\begin{pmatrix}
I									& (\mbA _{\Omega \Omega}^k)^{-1} \mbA_{\Omega \Gamma}^{k} \\
0									& I
\end{pmatrix}
\begin{pmatrix}
\ub_{\Omega}^{h,k} \\
\ub^{h,k}_\Gamma
\end{pmatrix}.
\end{align}
If we plug \eqref{eq:splittingBasisChange} into \eqref{eq:splittingBlockCholesky} 
we obtain
\begin{align}
    \langle \mbA^k  \begin{pmatrix} \ub_\Omega^{h,k}\\ \ub_\Gamma^{h,k} \end{pmatrix} , 
                    \begin{pmatrix} \ub_\Omega^{h,k}\\ \ub_\Gamma^{h,k} \end{pmatrix} \rangle = 
   \langle \mbA _{\Omega \Omega}^k   \vb_\Omega^{h,k} , 
                     \vb_\Omega^{h,k}\rangle 
    +
   \langle \mbS _k  \ub_\Gamma^{h,k}  , 
                    \ub_\Gamma^{h,k} \rangle  ,
\end{align}
where $\langle \cdot, \cdot \rangle$ is the standard Euclidean scalar product.
Similarly the forcing term can be transformed to
\begin{align}
    \langle  \begin{pmatrix} \tilde\fb_\Omega^{h,k}\\ \tilde\fb_\Gamma^{h,k} \end{pmatrix} , 
                    \begin{pmatrix} \ub_\Omega^{h,k}\\ \ub_\Gamma^{h,k} \end{pmatrix} \rangle = 
   \langle \tilde\fb^{h,k}_{\Omega}  , 
                     \vb_\Omega^{h,k}\rangle 
    +
   \langle \fb^h_k   , 
                    \ub_\Gamma^{h,k} \rangle ,
\end{align}
where
\begin{align}
\label{eq:subForcing}
&\fb^h_k =  \tilde{\fb}^{h,k}_{\Gamma} - \mbA_{\Gamma \Omega}^{k} (\mbA_{\Omega \Omega}^{k})^{-1} \tilde{\fb}^{h,k}_{\Omega}.
\end{align}
Hence, the energy functional $\mcE_k(u_k)$ can be written as
\begin{align*}
\mcE_k(u_k) =
\frac{1}{2} \langle \mbA_{\Omega \Omega}^k \vb^{h,k}_\Omega,\vb^{h,k}_\Omega \rangle
- \langle \tilde{\fb}^{h,k}_{\Omega}, \vb^{h,k}_\Omega\rangle
+ 	\underbrace
{
	\frac{1}{2} 
	\langle  \ub^{h,k}_\Gamma, \mbS^k  \ub^{h,k}_\Gamma \rangle
	- \langle \fb^h_k, \ub^{h,k}_\Gamma\rangle
}_{=:\mcE_k^\Gamma(\ub^{h,k}_{\Gamma})},
\end{align*}
We find that $\mcE_k^\Gamma(\cdot)$ does not depend on coefficients \(\ub_{\Omega}^{h,k}\) in $\Omega_k$.
So
the splitting
allows to find the nonlocal trace of the solution on $\Gamma$ by a constrained problem with the objective 
\begin{align}
\label{eq:SchurComplement}
\sum_{k=1}^K \mcE_k^\Gamma(\ub_\Gamma^{h,k}) = 
\frac{1}{2}\langle \mbS \ub^h, \ub^h \rangle - \langle \fb^h, \ub^h \rangle,
\end{align}
where $\mbS  = \diag(\family{\mbS^{k}})$,
$\ub^h = \family{\ub_\Gamma^{h,k}} \in U$,
and $\fb^h = \family{\fb^h_k} \in U$.
More precisely, the reduced, discrete formulation of \eqref{eq:subWeakNonlocal} is given by
\begin{equation}\label{eq:discreteSubNonlocal}
\begin{aligned}
&\mbox{given $\fb^h \in U$ find $\ub^h \in U$ such that} \\
&  \qquad
\ub^h = \argmin \limits_{\vb^h\in U} 
\frac{1}{2}\langle \mbS \vb^h, \vb^h \rangle - \langle \fb^h, \vb^h \rangle \\
&
\mbox{subject to $\Bb \ub^h = 0$}.
\end{aligned}
\end{equation}
There are no constraints imposed on nodes in $\Omega_k$ in problem \eqref{eq:subWeakNonlocal}.
So we can ultimately derive $\ub_\Omega^{h,k} \in U_{\Omega,k}$ from the backward substitution
\begin{align}
\label{eq:backwardSubstitution}
\ub_\Omega^{h,k} = (\mbA_{\Omega \Omega}^{k})^{-1} (\tilde{\fb}^{h,k}_{\Omega} - \mbA_{\Omega \Gamma}^k \ub_\Gamma^{h,k}  ),
\end{align}
which corresponds to the solution of $K$ independent Dirichlet problems.
The last step \eqref{eq:backwardSubstitution} requires no communication among the ranks. We therefore put it aside and focus only on
\eqref{eq:discreteSubNonlocal} in the following discussion.

\section{FETI}
\label{sec:FETI}
In the following, we briefly present a FETI algorithm which is commonly referred to as one-level FETI \cite{toselli2006domain, pechstein_finite_2013}. 
The algorithm has been developed for the solution of elliptic partial differential equations but it does not need to be changed for the nonlocal setting and the presentation works analogous to the classical, local case.

Problem \eqref{eq:discreteSubNonlocal} is a quadratic optimization problem with linear constraints.
The block diagonal matrix $\mbS$ is positive semi-definite and positive definite on $\ker(\Bb)$, i.e. $\ker(\mbS) \cap \ker(\Bb) = \lbrace 0 \rbrace$.
The necessary conditions of optimality therefore lead to the saddle point formulation
\begin{align}
\label{eq:discreteKKT}
\notag 
&\mbox{find $(\ub^h, \lambda) \in U \times \Lambda$ such that} \\
&\begin{cases}
& \mbS \ub^h + \Bb ^\top \lambda = \fb^h \\
& \Bb       \ub^h	  = 0.
\end{cases}
\end{align}
Note that $\ker (\Bb^\top) = \lbrace 0 \rbrace$ and $\range(\Bb) = \Lambda$ as we assume that $\Bb$ has full row-rank.
The block diagonal matrix $\mbS$ contains singular blocks $\mbS_k$ whenever $\widehat{\Omega}_k$ is a floating subdomain.
For scalar kernels the block is 1-dimensional and consists of constant vectors.
In the case of 2d tensorial kernels it is 3-dimensional.
Let now $\NFloating$
be the number of floating subdomains and let
$\Zb_k$ be such that $\range(\Zb_k) = \ker(\mbS_k)$. 
We then construct a full-column-rank matrix
\begin{align}
\label{eq:matKerS}
\Zb = \begin{pmatrix}
\Zb_1 	& 0 		& \cdots & 0	 \\
0 		& \Zb_2 	& \ddots & \vdots \\
\vdots 	& \ddots	& \ddots & 0	 \\
0 		& \cdots	& 0		 & \Zb_{\NFloating} \\
\end{pmatrix},
\end{align} 
so that $\range (\Zb) = \ker(\mbS)$. The exact knowledge of the kernel of $\mbS$ allows an efficient evaluation of the pseudoinverse.

The system in \eqref{eq:discreteKKT} can be reduced to a problem depending only on $\lambda$. To that end we derive two equations. 
First of all we find that the first line
in \eqref{eq:discreteKKT} is solvable only if $\fb^h - \Bb^\top\lambda \in \range(\mbS) = \ker(\Zb^\top)$. Hence,
$\Zb^\top (\fb^h - \Bb^\top \lambda) = 0$ which is equivalent to
$
\Gb^\top \lambda = \Zb^\top \fb^h,
$
where $\Bb \Zb$ is denoted by $\Gb$\label{inline:G}. Secondly, by the help of a pseudoinverse, and the fact that $\range(\mbS)^\perp = \range(\Zb)$,  the first equation in \eqref{eq:discreteKKT} can be written as
\begin{align}
\label{eq:forU}
\ub^h = \mbS^+(\fb^h - \Bb^\top \lambda) - \Zb \alpha,
\end{align}
for some unknown vector $\alpha$. 
We choose here the Moore-Penrose generalized inverse so that $\mbS^+$ is symmetric.
If we plug in that equation into the second line of \eqref{eq:discreteKKT}
we get the system 
\begin{align}
\label{eq:FETIii}
\Bb \mbS^+ \Bb^\top \lambda + \Gb \alpha = \Bb \mbS^+\fb^h. 
\end{align}
Equation \eqref{eq:FETIii} depends not only on $\lambda$ but also on an unknown $\alpha$. However, we can again separate those problems. To that end we eliminate the term $\Gb\alpha$ in \eqref{eq:FETIii} with the help of the 
orthogonal projection
\begin{align}
\label{eq:P}
\Pb = \Ib - \Gb(\Gb^\top \Gb)^{-1} \Gb^\top
\end{align}
and determine the unknown $\alpha$ for a given $\lambda$ in a second step by
\begin{align}
\label{eq:forAlpha}
\alpha = (\Gb^\top \Gb)^{-1 }\Gb^\top (\db - \Fb \lambda).
\end{align}
Finally, we obtain the following system of equations 
\begin{align}
\label{eq:FETIproblem}
\begin{cases}
& \Pb \Fb \lambda = \Pb \db \\
& \Gb ^\top \lambda = \eb,
\end{cases}
\end{align}
where $\Fb = \Bb \mbS^+ \Bb^\top$, $\db = \Bb \mbS^+\fb^h$, and $\eb = \Zb^\top \fb^h$. 
One-level FETI is a distributed preconditioned projected conjugate gradient method to solve \eqref{eq:FETIproblem}.
Clearly, the evaluation of the various matrices in the system requires special care and we postpone the discussion of their numerical implementation to Section \ref{sec:implementation}.
The increments of the projected cg-method lie in $\range(\Pb) = \ker(\Gb^\top) = \Lambda_{ad}$ \label{inline:admissibleIncrements}
and we see that any $\tilde{\lambda} \in \Lambda_{ad}$ fulfills $\Gb^\top \tilde{\lambda} = 0$. 
Hence, it suffices to choose a starting value $\lambda^0$ with $\Gb^\top \lambda^0  = \eb$ and the 
solution of the projected
cg-algorithm ultimately solves \eqref{eq:FETIproblem}. 
In Algorithm \ref{alg:oneLevelFETI} we summarize a preconditioned version of one-level FETI.
\begin{algorithm}
	\caption{One-level FETI} 
	\label{alg:oneLevelFETI}
	\begin{enumerate}
		\item build substructures and construct $\zeta$, $\Zb$, $\Bb$, and $\Bb_D$.
		\label{alg:oneLevelFETI:substructures}
		\item assemble in parallel $\mbA^k$ and $\fb^h_k$
		\label{alg:oneLevelFETI:assemble}
		\item compute $\db$, $\eb$ and $\Gb$
		\item compute Cholesky decompositions to evaluate $\mbS^+$ and $(\Gb^\top \Gb)^{-1}$
		\label{alg:oneLevelFETI:cholesky}
		\item set up $\Fb$, the projection $\Pb$, and the preconditioner $\Mb^{-1}$
		\label{alg:oneLevelFETI:preconditioner}
		\item initialize $\lambda^0 = \Gb(\Gb^\top \Gb)^{-1} \eb$, \eqref{eq:FETIproblem}
		\label{alg:oneLevelFETI:lambda}
		\item solve $\Pb \Fb \lambda = \Pb \db$ with projected pcg and initial value $\lambda^0$,
		\label{alg:oneLevelFETI:pcg} 
		\eqref{eq:FETIproblem}
		\item solve $\alpha = (\Gb^\top \Gb)^{-1 }\Gb^\top (\db - \Fb \lambda)$,
		\label{alg:oneLevelFETI:alpha} 
		\eqref{eq:forAlpha}
		\item solve $\ub^h = \mbS^+(\fb^h - \Bb^\top \lambda) - \Zb \alpha$, \eqref{eq:forU}
		\label{alg:oneLevelFETI:u} 
	\end{enumerate}
\end{algorithm}
\subsection{Preconditioning}
The preconditioner we used in the numerical tests is given by
\begin{align}
\label{eq:preconditioner}
\Mb^{-1} = \Bb_D ~ \mbS ~  \Bb_D ^\top,
\end{align}
where 
$\Bb_D = (\Bb  \Db^{-1} \Bb ^\top )^{-1}  \Bb \Db^{-1} $.
Here each block $\Db^k$ of the diagonal matrix $\Db$ counts the number of times a node $\xb^h_m \in \Gamma_k^{N,h}$ occurs, i.e.
$(\Db^k)_{m m} = \zeta(\xb^h_m, \xb^h_m)$.
Commonly, $\Mb^{-1}$ denoted as \textit{scaled Dirichlet preconditioner} \cite{toselli2006domain, pechstein_finite_2013}.

\section{Implementation}
\label{sec:implementation}
In our implementation we use \mbox{PETSc} \cite{balay2019petsc} which at its core provides MPI distributed linear algebra.
Therefore we only need to describe how the spaces $U$, and $\Lambda$ are distributed among MPI ranks. This immediately
defines the layout of all vectors in the spaces. A matrix is distributed according to the layout of its range space. 
Our implementation guarantees that the parallel layout of all vectors in $U$ and all matrices mapping to $U$ is aligned 
according to the decomposition of domains. 
That is, each $U_k$ is assigned to one rank so that the evaluation of the block diagonal operators $\mbS$ and $\mbS^+$ requires no communication.
Consequently, the number of subdomains $K$ defines the number of MPI ranks.
Furthermore, \mbox{PETSc} automatically distributes vectors and matrices belonging to the spaces $\Lambda$.
The communication between the processes, e.g. within the projected conjugate gradient method, is also automatically handled by \mbox{PETSc} for the given data layout.

We now describe the assembly and evaluation of the matrices, and the settings of the applied solvers. To this end 
Section \ref{subsec:buildingSubtructures} to Section \ref{subsec:projectedCgMethod} follow the most important steps in Algorithm \ref{alg:oneLevelFETI}. In Table \ref{tab:implementation} we give a summary of the implemented matrices and functions.

\subsection{Building substructures}
\label{subsec:buildingSubtructures}

The construction of substructures, see Algorithm \ref{alg:oneLevelFETI} Line \ref{alg:oneLevelFETI:substructures}, is based on an estimate of the support of the kernel so that the covering fulfills \eqref{eq:subRequirement} and does not cut through finite elements.
In order to simplify the task of constructing substructures, we will be using approximate interaction balls such that \(B_{\delta,p}^{h}(\xb)\subset B_{\delta,p}(\xb)\) as it allows to directly work with the \(\ell^{p}\)-norm to decide whether two points might interact via the kernel \(\gamma_{h}\).
Following \cite{capodaglio2020general} we can concurrently construct the nonlocal subdomains $\widehat{\Omega}_k$ based on a local decomposition 
$\family{\widetilde{\Omega}_k }$ of the mesh.
The local decomposition can be obtained from a graph partitioning program like METIS \cite{metis}. In our experiments we instead subdivide the domain into rectangles which is easier to implement but obviously not as general. However, it allows for similar experimental settings for meshes of different resolution and gives easy control over the number of floating subdomains.
In general, given $\family{\widetilde{\Omega}_k }$, we can construct the nonlocal subdivision $\family{\widehat{\Omega}_k}$ by adding an approximate $\delta/2$-\textit{neighborhood} to $\widetilde{\Omega}_k$ \cite{capodaglio2020general}. If the approximate neighborhood covers the interaction set of the kernel it can be shown that the resulting extensions $\widehat{\Omega}_k$ yield a nonlocal subdivision \cite{capodaglio2020general}, that is that the multidomain formulation is equivalent to the original problem.

We give a short description of the algorithm which extends a given local subdivision to a nonlocal one.
The procedure determines the sets $\widehat{\Omega}_k$ concurrently on $K$ parallel ranks. Let us fix rank $k$ where the nonoverlapping subdomain $\widetilde{\Omega}_k$ is stored.
At first the algorithm identifies the elements $E^{k'} \subset\widetilde{\Omega}_k$ which have a non empty intersection with the boundary $\partial\widetilde{\Omega}_k \cap \partial\widetilde{\Omega}_{k'}$ for each neighboring subdomains \(k'\).
Using this initial set of elements in a \(\delta/2\)-neighborhood of \(\widetilde{\Omega}_{k'}\) the algorithm greedily checks and adds elements in \(\widetilde{\Omega}_k\) in order of their distance in the mesh connectivity graph.
The search terminates when no more elements in the \(\delta/2\)-neighborhood of \(\widetilde{\Omega}_{k'}\) can be found, and the element information is communicated to subdomain \(k'\).
After receiving mesh elements from all its neighbors the nonlocal domain \(\widehat{\Omega}_k\) is constructed by concatenation.

Again, following \cite{capodaglio2020general}, the subsets $\Gamma_k^D$ of $\Gamma^D$ 
are defined as interaction domain of the sets $\widetilde{\Omega}_k$. Hence, $\Gamma_k^D$ contains all elements in $\Gamma^D$ which lie in a $\delta$-neighborhood of some element in $\widetilde{\Omega}_k$. This results in larger overlaps among adjacent sets $\Gamma_k^D$ and $\Gamma_k^D$, see the darker blue subregion of $\Gamma_1^D$ in Figure \ref{fig:subdomains}.

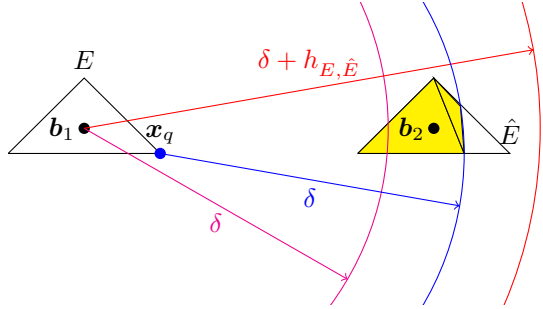
\begin{figure}
	\centering
	\begin{tikzpicture}[scale=2]
		
		\def\deltaVal{2.}
		\def\height{0.5}
		\def\hVal{1.}
		\def\dist{1.3}
		
		\coordinate (bary1) at (-0.5*\dist-0.5*\hVal,0.333*\height);
		\coordinate (bary2) at (0.5*\dist+0.5*\hVal,0.333*\height);
		
		\clip (-\dist-\hVal,-2*\height) rectangle + (2*\dist+2*\hVal,4*\height);
		
		\coordinate (elem11) at (-0.5*\dist-\hVal,0);
		\coordinate (elem12) at (-0.5*\dist,0);
		\coordinate (elem13) at (-0.5*\dist-0.5*\hVal,\height);
		\draw (elem11) -- (elem12) -- (elem13) -- (elem11);
		\coordinate (elem21) at (0.5*\dist,0);
		\coordinate (elem22) at (0.5*\dist+\hVal,0);
		\coordinate (elem23) at (0.5*\dist+0.5*\hVal,\height);
		\draw (elem21) -- (elem22) -- (elem23) -- (elem21);

		\draw[fill=yellow] (elem21) -- ($(elem12) + (\deltaVal,0)$) -- (elem23) -- (elem21);
		\draw[fill=yellow] ($(elem12) + (\deltaVal,0)$) -- ($(elem12) + (9:\deltaVal)$) -- (elem23) -- ($(elem12) + (\deltaVal,0)$);
		
		\node at (bary1) [circle,fill,inner sep=1.5pt] {};
		\node at (bary1) [anchor=east] {\(\bm{b}_{1}\)};
		\draw[color=magenta] (bary1) circle (\deltaVal);
		\draw[->,color=magenta] (bary1) -- node[below] {\(\delta\)} ($ (bary1) + (-30:\deltaVal) $);
		
		\draw[color=red] (bary1) circle (\deltaVal+\hVal);
		\draw[->,color=red] (bary1) -- node[above] {\(\delta+h_{E,\hat{E}}\)} ($ (bary1) + (10:\deltaVal+\hVal) $);
		
		\node[color=blue] at (elem12) [circle,fill,inner sep=1.5pt]{};
		\node at (elem12) [anchor=south] {\(\xb_{q}\)};
		\draw[color=blue] (elem12) circle (\deltaVal);
		\draw[->,color=blue] (elem12) -- node[below] {\(\delta\)} ($ (elem12) + (-10:\deltaVal) $) ;
		
		\node at (bary2) [circle,fill,inner sep=1.5pt]{};
		\node at (bary2) [anchor=east] {\(\bm{b}_{2}\)};
		
		\node at (elem13) [anchor=south] {\(E\)};
		\node at (elem22) [anchor=south] {\(\hat{E}\)};
		
	\end{tikzpicture}
	\caption{
		Left: An element \(E\) of the outer integration; right: an element \(\hat{E}\) of the inner integral. \\
		In {\color{blue}blue}: \(\delta\)-ball around a quadrature node \(\xb_{q}\) of the outer integration.\\
		In {\color{yellow}yellow}: two sub-simplices of the inner element that are part of the approximate ball around \(\xb_{q}\).\\
		In {\color{red}red} and {\color{magenta}magenta}: \((\delta+h_{E,\hat{E}})\)- and \(\delta\)-balls around the barycenter \(\bm{b}_{1}\) of the left element.
		It can be seen that \(\delta<|\bm{b}_{1}-\bm{b}_{2}|<\delta+h_{E,\hat{E}}\).
	}
	\label{fig:barycenterCriterion}
\end{figure}

The exact criterion for deciding whether an element \(E\) is in the \(\delta/2\)-neighborhood of \(\widetilde{\Omega}_{k'}\) is \(\inf_{\xb\in E,\yb\in\widetilde{\Omega}_{k'}} ||\xb-\yb||_{2} \leq \delta/2\).
This can be relaxed to computationally less expensive variants:
in \cite{capodaglio2020general} the authors propose to consider the $\ell_2$-ball around the element barycenters of radius $\delta/2 + h$.
The barycenter criterion can be checked easily but it requires to add the mesh size parameter $h$ which might be unsuitable for nonuniform grids. However, nonuniform grids can be handled by replacing $h$ by the element dependent value
$h_{E, \widehat{E}}$ 
equal to the maximal diameter of $E^{k'}$ and $\widehat{E}^{k'}$.
We illustrate the barycenter criterion in Figure~\ref{fig:barycenterCriterion}.

Alternatively, one can also test whether the minimum distance of element vertices is less than \(\delta/2\) to determine interaction of elements.

\begin{figure}
	\begin{subfigure}{0.49\linewidth}
		\includegraphics[page=1, width=1.\linewidth]{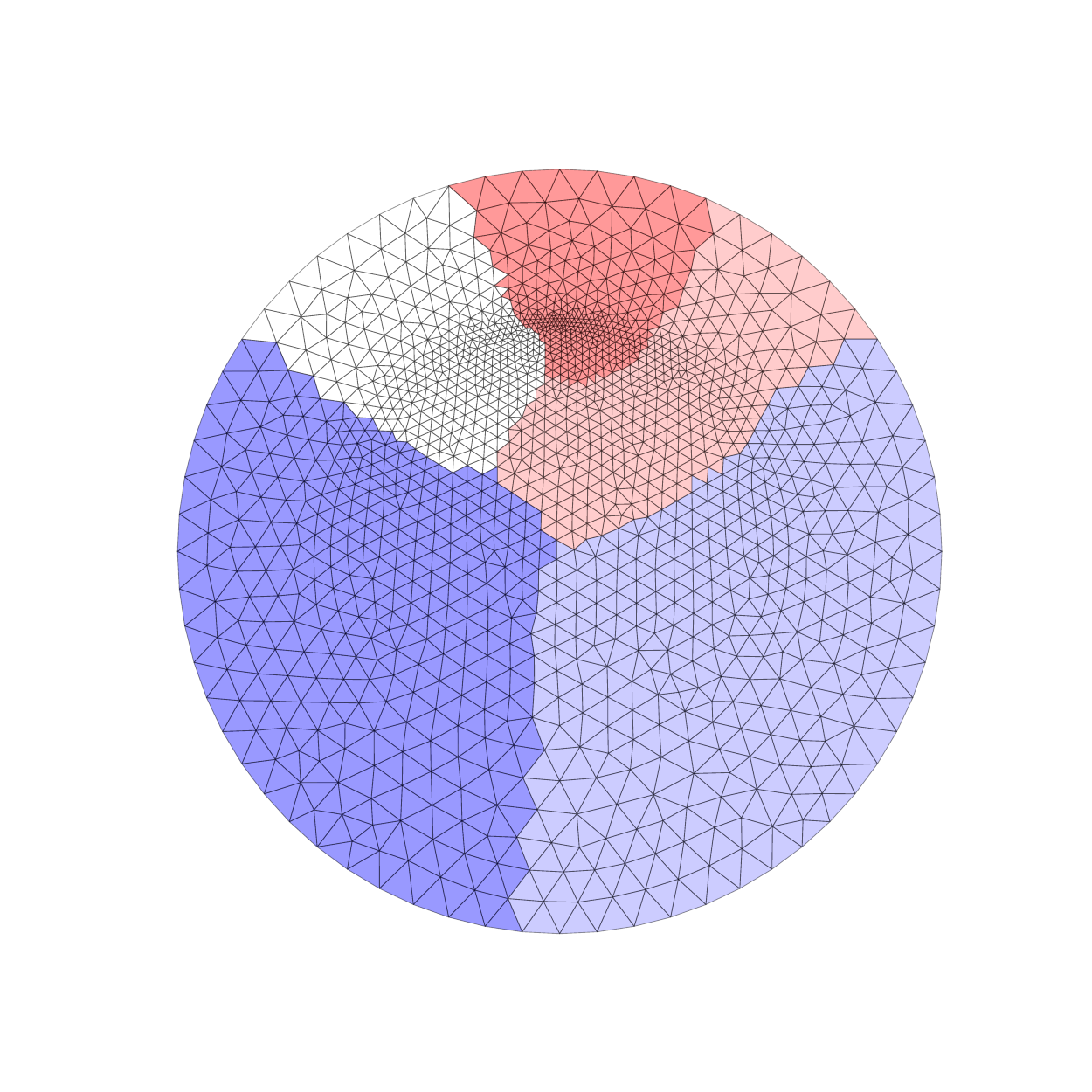}
		\caption{}
		\label{fig:nlsubdivision:a}
	\end{subfigure}
	\begin{subfigure}{0.49\linewidth}
		\includegraphics[page=2, width=1.\linewidth]{Figures/nlsubdivision}
		\caption{}
		\label{fig:nlsubdivision:b}
	\end{subfigure}
	\caption{(a) Decomposition of a nonuniform mesh into 5 nonoverlapping subdomains \(\family{\widetilde{\Omega}_{k}}\).
          (b) The counting function \(\zeta\) corresponding to a nonlocal subdivision $\family{\widehat{\Omega}_k}$.
          The white and red elements mark the interaction set $\Gamma$.
          The blue elements belong to one subdomain, the white elements to two, and the red elements to three.}
	\label{fig:nlsubdivision}		
\end{figure}

Figure \ref{fig:nlsubdivision:a} shows a nonuniform mesh inside a disk which is partitioned into 5 subdomains $\widetilde{\Omega}_k$.
Figure \ref{fig:nlsubdivision:b} shows the overlaps of the corresponding nonlocal subdivision where the alternative $\delta/2$-neighborhoods has been applied.

Based on the decomposition $\family{\widehat{\Omega}_{k}\cup\Gamma^D_{k}}$ of we define $\zeta$, $\Zb$, $\Bb$, and $\Db$.
The function $\zeta$ can be evaluated at vertices $\xb_\ell, \xb_{m} \in \Omega \cup \Gamma^D$ without much effort. Note that it might exhibit jumps at some vertices, but is always well-defined, see Remark \ref{rem:zetajumps}.
To that end, we store the values of $\ind_{\Omega_k\cup\Gamma^D_k}(\xb_\ell)$ for 
$\xb_\ell \in \Omega \cup \Gamma^D$
in a matrix $(\Cb)_{\ell k}$ and obtain, by definition, that $\zeta(\xb_\ell, \xb_{m}) = (\Cb \Cb^\top)_{\ell m}$.
The number of nonzero entries in $\Cb$ grows only mildly with increasing number of subdomains for a fixed mesh size, so that it is sparse.
A single entry $(\Cb \Cb^\top)_{\ell m}$ of the sparse matrix-matrix product can be computed at low cost as the number of overlaps $\max_{\ell, m} (\zeta(\xb_\ell, \xb_m))$  is bounded independent of $K$.
The counting function $\zeta$ is constant on (the interior of) finite elements as the decomposition does not cut through elements. 
Therefore, an analogous procedure can be used to evaluate $\zeta$ for any pair of elements. This is needed for the assembly of the matrices $\mbA^k$.
The matrix $\Zb$ contains the rigid body modes of the floating subdomains. For scalar kernels it is immediately obtained from $\Cb$. In the tensorial case we additionally need to provide a vector with coefficients 
corresponding to a rotation of each floating subdomain.

Due to \eqref{eq:constraintSpace} we can assess the number of constraints $M_C$ in the non-redundant case with the help of $\zeta$. 
Furthermore each row of $\Cb$ allows to identify all local indices $m_k$ for $k=1,\dots,K$ corresponding to the global index $m$ of an interface node $\xb_m \in \Gamma^h$
and we know that each row of $\Bb$  contains exactly two entries, $1$ and $-1$.
The construction of $\Bb$ in a sparse format is then straightforward.
The diagonal matrix $\Db$ is defined based on $\zeta$, see \eqref{eq:preconditioner}. 
Finally, we compute and store the scaled matrix of constraints $\Bb_D  = (\Bb  \Db^{-1} \Bb ^\top )^{-1}  \Bb \Db^{-1}$. 
Note that $\Bb  \Db^{-1} \Bb ^\top $ is block diagonal where the block size is smaller than $\max_{\ell} (\zeta(\xb_\ell, \xb_\ell)) \ll K$.

\subsection{Assembly of the stiffness matrix} 
Given the decomposition we distribute the relevant parts of $\Cb$, that way each rank can compute the value
of $\zeta$ on its subdomain and the local stiffness matrices
$\mbA^k$ can be assembled concurrently, see Line \ref{alg:oneLevelFETI:assemble} of Algorithm \ref{alg:oneLevelFETI}.
However, an efficient assembly which allows large scale experiments is a challenge on its own.
We discuss selected aspects concerning the three kernel functions which are used in the numerical experiments.
More details about the implementation of nonlocal operators with truncated interaction horizon can be found e.g. in \cite{nlfem}.
In general kernel $\ell_\infty$-balls truncations can be evaluated without geometric error.
For $\ell_2$-ball truncations we apply an error commensurate 
polygonal approximation of the support \cite{DEliaFEM2020}.
The assembly the operators based on symmetric kernels leads to symmetric stiffness matrices and the null spaces of floating subdomains
contain the rigid body modes up to machine precision. 

\subsection{Pseudoinverse and Cholesky decomposition}
In Line \ref{alg:oneLevelFETI:cholesky} of Algorithm \ref{alg:oneLevelFETI}
we prepare the evaluation of $\mbS^+$ and $(\Gb^\top \Gb)$.
To that end, consider $\mbA^{k}$ as given in \eqref{eq:subStiffnessMatrix} and let $\vb^h_k = (0, \vb^{h,k}_\Gamma) \in U_{\Omega, k} \times U_{\Gamma, k}$ be some vector of coefficients.
The evaluation $\mbS^+_k \vb^{h,k}_\Gamma = \wb^{h,k}_{\Gamma}$ is obtained from a solution of the Neumann problem
\begin{align}
\label{eq:neumannSubproblem}
\begin{pmatrix}
\mbA_{\Omega \Omega}^{k} 		& \mbA_{\Omega \Gamma}^{k} \\
\mbA_{\Gamma \Omega}^{k}	& \mbA_{\Gamma \Gamma}^{k}
\end{pmatrix}
\begin{pmatrix} 
\wb_{\Omega}^{h,k} \\ 
\wb^{h,k}_{\Gamma} \end{pmatrix}
= \begin{pmatrix}
0 \\  
\vb^{h,k}_\Gamma
\end{pmatrix},
\end{align}
with $\wb_{\Omega}^{h,k} \in U_{\Omega,k}$ and $\wb_{\Gamma}^{h,k} \in U_{\Gamma,k}$.
If the domain $\Gamma^D_k$ has positive measure this system has a unique solution for any $\vb^{h,k}_\Gamma \in U_k$. 
We can then immediately solve the system by means of a Cholesky decomposition. 
In our implementation we use \mbox{CHOLMOD} \cite{chen2008algorithm}.
If otherwise 
$\widehat{\Omega}_k$ is a floating subdomain we evaluate its pseudoinverse, where we can assume that problem \eqref{eq:neumannSubproblem} has a solution.
In order to use a direct solver which does not accept singular matrices, artificial zero Dirichlet boundary conditions are imposed.
We explain the procedure in more detail for the case of the scalar nonlocal diffusion problem.
In this case, the nullspace of \(\mbA^{k}\) for floating subdomains is given by the constant function.
We choose some node $\ell$ and replace the original system by $\widetilde{\mbA}^{k}$ where the $\ell$-th row and column are replaced by the $\ell$-th standard unit vector.
We now set the $\ell$-th element of the right hand side $\vb^h_k$ to zero and denote the resulting vector by $\widetilde{\vb}^{h}_k$.
This amounts to enforcing a homogeneous Dirichlet condition in a single point, thereby making the matrix $\widetilde{\mbA}^{k}$ invertible,
hence we can find
$
\widetilde{\mbA}^{k}
\widetilde{\wb}^{h}_k
= 
\widetilde{\vb}^{h}_k.
$
The matrices $\mbA^{k}$ and $\widetilde{\mbA}^{k}$ differ in the $\ell$-th column but $(\wb^h_k)_\ell = 0$, hence it is clear that $(\mbA^{k} \widetilde{\wb}^{h}_k)_m = (\vb^h_k)_m$ for all $m\neq \ell$.
As the solution using the Moore-Penrose generalized inverse $\vb^h_k$ has zero mean, this yields that
\begin{align*}
(\vb^h_k)_\ell & 
= - \sum_{m\neq \ell} (\vb^h_k)_m 
= - \sum_{m\neq \ell} (\widetilde{\vb}^{h}_k)_{m}
= - \sum_{m\neq \ell} (\widetilde{\mbA}^{k} \widetilde{\wb}^{h}_k)_m = \\
& = - \sum_{m\neq \ell} \sum_{i} (\mbA^{k})_{m,i} (\widetilde{\wb}^{h}_k)_i
= - \sum_{i} (\widetilde{\wb}^{h}_k)_i \sum_{m\neq \ell} (\mbA^{k})_{m,i} \\
& = \sum_{i} (\mbA^{k})_{\ell,i} (\widetilde{\wb}^{h}_k)_i
  = (\mbA^{k} \widetilde{\wb}^{h}_k)_\ell.\end{align*}
In the second to last step, we have used that the constant function is in the nullspace of \(\mbA^{k}\) and that \(\mbA^{k}\) is symmetric.
Hence, $\widetilde{\wb}^{h}_k$ solves \eqref{eq:neumannSubproblem}. We finally subtract the mean from $\widetilde{\wb}^{h}_k$ to obtain a solution which is orthogonal to $\ker(\mbS)$,
so that we have evaluated the Moore-Penrose generalized inverse 
$\mbS^+ \vb^{h,k}_\Gamma$. 
Note that the overall procedure for singular matrices is not more expensive than a regular solve, apart from some vector operations.
An analogous procedure is used for the case of tensorial kernels.

The projection $\Pb$ involves an evaluation of $(\Gb^\top \Gb)^{-1}$. For scalar kernels the system size is 
$\NFloating \times \NFloating$, and $3 \NFloating \times 3 \NFloating$ for tensorial kernels, so comparably small.
The matrix $\Gb^\top \Gb$ is furthermore symmetric and sparse. 
We use \mbox{CHOLMOD} \cite{chen2008algorithm} for the evaluation of its inverse. Furthermore, the size of the system allows to redundantly solve it on all ranks which saves communication.

\subsection{Preconditioned solver}
\label{subsec:projectedCgMethod}
The projected CG-method 
to solve $\Pb \Fb \lambda = \Pb \db$ is performed by the PETSc-function \texttt{KSPSolve()}.
The preconditioner requires an evaluation of the block diagonal matrix $\mbS$. One evaluation of $\mbS_k$ requires a solve of the local Dirichlet problem \eqref{eq:subSchurComplement} on each subdomain. The matrices 
$\mbA_{\Omega \Omega}^k$ are invertible, irrespective of the regularity of $\mbA^k$. Therefore $k$ concurrent Cholesky decompositions, again using \mbox{CHOLMOD} \cite{chen2008algorithm},  allow a direct evaluation of $\mbS$.
The preconditioner $\Mb^{-1}$ furthermore requires the application of $\Bb_D$ and $\Bb_D^\top$ which have both been set up in the first step of Algorithm \ref{alg:oneLevelFETI}.

\begin{table}
	\begin{tabular}{lll}
		\textbf{Name} 											& \textbf{Format}	& \textbf{Details} 													\\
		$\Cb, \Zb, \Bb, \Bb_D$									& stored			& sparse format														\\
		$\zeta(\xb_\ell, \xb_{\ell'})$							& evaluated			& $(\Cb \Cb^\top)_{\ell \ell'}$										\\
		$\mbA_k, \mbA_{\Omega\Omega}^k, \dots$		& stored 			& sparse format														\\
		$\mbA_k^+$											& factorization stored		& CHOLMOD, possibly rank-1 correction								\\
		$\mbS^+_k$											& evaluated 		& solve \eqref{eq:neumannSubproblem}, requires $\mbA_k^+$		\\
		$\Fb$													& evaluated			& requires $\Bb$, $\mbS^+$, and 
		$\Bb^\top$										\\
		$\Gb, \Gb^\top$											& stored & sparse format														\\
		$\Gb^\top \Gb$											& stored redundantly& sparse format														\\
		$(\Gb^\top \Gb)^{-1}$									& fact. stored	redundantly		& CHOLMOD													\\
		$\Pb$													& evaluated redundantly	& requires $\Gb, \Gb^\top$ and $(\Gb^\top \Gb)^{-1}$				\\
		$(\mbA_{\Omega\Omega})^{-1}$							& factorization stored			& CHOLMOD															\\
		$\mbS_k$											& evaluated 		& evaluate \eqref{eq:subSchurComplement}, requires 
		$(\mbA_{\Omega\Omega})^{-1}$									\\	
		$\Mb^{-1}$									& evaluated			& requires  $\Pb$, $\Bb_D$, and $\mbS_k$						\\		
	\end{tabular}
	\caption{Numerical implementation of functions and matrices.}
	\label{tab:implementation}
\end{table}

\section{Numerical experiments}
\label{sec:numerical-experiments}
\subsection{Kernels used for the numerical illustrations}\label{sec:kernels}
For the numerical illustrations,
we restrict our attention to the two-dimensional case $d=2$. The kernels can be scalar or matrix-valued.

For the nonlocal diffusion setting, we consider two kernels. First, we have the scalar-valued kernel that is {\em constant} within its support given by
\begin{equation}
	\label{eq:constantkernel}
	\gamma_{c}(\xb,\yb) = \phi_{c}(\xb,\yb)  
	\ind_{B_{\delta,\infty}(\xb)}(\yb)
	\,\,\,\text{with}\,\,\, \phi_{c}(\xb,\yb) = \frac{3}{4 \delta^4}.
\end{equation}
The scaling constant ${3}/{4 \delta^4}$ is chosen so that that the nonlocal diffusion operator, i.e., the left-hand side of (\ref{eq:strong_nonloc}a), converges as $\delta\to0$ to the (negative of the) Laplace operator \cite{DEliaFEM2020}. We also consider, for $s\in (0,1)$, the scalar-valued kernel that is of {\em fractional-type} within its support given by
\begin{equation}
	\label{eq:fractionalkernel}
	\gamma_s(\xb, \yb)  = \phi_{s}(\xb,\yb)
	\ind_{B_{\delta,2}(\xb)}(\yb)	
	\,\,\,\text{with}\,\,\, \phi_{s}(\xb,\yb)	=
	\frac{2-2s}{\pi \delta^{2-2s}}
	\frac{1}{|\xb-\yb|^{d + 2s}_2}.
\end{equation} 
Again, the scaling constant $(2-2s)/{\pi \delta^{2-2s}}$ is chosen so that the nonlocal diffusion operator converges as $\delta\to0$ to the (negative of the) Laplace operator.  
We refer to \eqref{eq:fractionalkernel} as a ``fractional-type kernel'' because the singularity of that kernel is the same as that of the fractional Laplacian kernel which is positive on all of $\R^2$; see \cite{acta20}.

For the solid mechanics setting, we consider the matrix-valued nonlocal kernel
\begin{equation}
	\label{eq:peridynamics}
	\gammab_{PD}(\xb, \yb) =   \Phib_p(\xb,\yb)
	\ind_{B_{\delta,2}(\xb)}(\yb)
	\,\,\,\text{with}\,\,\,
	\Phib_p(\xb,\yb)	=   \frac{3}{\delta^3}
	\frac{(\xb - \yb) \otimes (\xb - \yb)}{|\xb - \yb|^3}
\end{equation}
which is a special case of the kernel which defines the bond-based {\em peridynamics} model for solid mechanics introduced in \cite{silling2000reformulation}. The constant $3/\delta^3$ is chosen so that that the nonlocal solid mechanics operator, i.e., the left-hand side of (\ref{eq:strong_nonloc-v}a), converges as $\delta\to0$ to the constant coefficient local Navier partial differential operator
$
-\frac{\pi}{4}\Delta \ub - \frac{\pi}{2} \nabla (\nabla\cdot\ub)
$ 
for a Poisson ratio $1/4$ \cite{mengesha2014bond}.


\subsection{Problem setting}
As domain we set $\widehat{\Omega} = (0,1)^2$. 
We impose no Neumann constraints so that $\Gamma_D$ coincides with the
interaction domain and $\Gamma^D = [-\delta, 1 + \delta]^2 \setminus \widehat{\Omega}$ for $\gamma_{c}$.
Note that the interaction domains of $\gamma_s$ and $\gammab_{PD}$ are also contained in $\Gamma_D$ as the $\ell_2-$ball is contained in the $\ell_\infty-$ball of equal radius.
For the diffusion problems we impose the forcing term $f(\xb) = -2(1 + \xb_1)$ and obtain the corresponding exact solution
$u(\xb) = \xb_1^2 \xb_2 + \xb_2 ^2$. 
The peridynamic problem is solved for the forcing term $\fb(\xb) = - \frac{\pi}{2}(1 + 2\xb_1 , \xb_2 )^\top$ 
and the exact solution $\ub(\xb) = (\xb_2^2, \xb_1^2 \xb_2)$.
We choose a regular triangulation in $\widehat{\Omega}$.
A regular mesh in $\widehat{\Omega}$ allows a precise determination of the rate of convergence to the true solution for vanishing mesh size $h$ but the FETI approach, and our implementation in particular, do not rely on it.
The domain $\widehat{\Omega}$ is substructured into $K = k_1 \times k_2$ rectangles. See Figure \ref{fig:mesh-decomposition} for an example of a subdivision into $3 \times 3$ subdomains. 

\subsection{Computational resources}
All experiments were performed on the Solo
cluster at Sandia National Laboratories \cite{sandiaHPC}. Solo consists of 374 nodes with
dual socket 2.1 GHz Intel Broadwell CPUs. Each node has 36 cores and 128
GB of RAM. Solo uses an Intel Omni-Path interconnect.

\subsection{Considered solver strategies}

In what follows, we present numerical results using two solvers.
For a performance baseline, we use a distributed conjugate gradient method with diagonal preconditioner to solve problem \eqref{eq:weak_nonloc}.
The rows of the discrete set of equations is distributed over MPI ranks using a partition that is similar to the subdomain distribution of the FETI solver, but none of the unknowns are duplicated.
This mimics the default solver approach used by the Peridigm project~\cite{peridigm}.
We compare this baseline approach with the FETI solver presented above both in terms of iterations and in terms of solution time.
We do not include the setup time of the solvers in the solution time.
In practice, setup of the FETI solver takes longer than a simple Krylov solver, but will easily be amortized if the system of equations needs to be solved repeatedly.
It will also become apparent that the conjugate gradient method is not always scalable, which means that for fine enough discretization FETI will outperform the Krylov method for a single solve even when setup cost is taken into account.

\subsection{Weak scaling for fixed horizon}
\label{sec:weak-scaling-fixed}

In a first example, we fix the horizon \(\delta=0.008\) and vary the mesh size \(h\).
The number of subdomains is chosen proportional to the global number of unknowns.
This also implies that the size of the subdomain matrices \(\mbA^k\) does not change.
The number of non-zero entries per row of \(\mbA^k\) is proportional to \((\delta/h)^{2}\) and hence increases as \(h\rightarrow0\).

In Table~\ref{tab:fixedHorizon} we display the resulting iteration counts and solve timings for the two solvers, as well as the computed \(L^{2}\) error that matches for both solutions.
In all cases, we set a convergence tolerance of \(10^{-10}\).
In the first case of the constant kernel \(\gamma_{c}\), we observe that both diagonally preconditioned CG and FETI solvers take a number of iterations that is essentially independent of the problem size, and that the solve time increases in both cases.
The increase in time is caused by the density of the operators, whereas the constant number of iterations taken by CG is explained by the fact that the condition number of the global matrix \(\mbA\) is constant \cite{aksoylu2011variational}.

We can contrast this behavior with the case of the fractional kernel \(\gamma_{s}\) with \(s=0.4\).
It is observed that the number of iterations for the Krylov method is no longer constant, but grows as the mesh is refined.
This is explained by the fact that the condition number of the matrix \(\mbA\) grows as \(h^{-2s}\) \cite{aksoylu_conditioning_2014}.
The advantage of the FETI approach becomes apparent, as the solution time displays more favorable scaling.

For the peridynamic kernel \(\gammab_{PD}\), CG is a scalable solver in terms of iterations, yet the benefit of the FETI approach in terms of timings is apparent.

\begin{table}
	\centering
	
	\begin{tabular}{lrr|rr|rr|rr}
		\multicolumn{9}{c}{Constant kernel \(\gamma_{c}\)} \\
		\hline
		&&&\multicolumn{2}{c|}{CG}&\multicolumn{2}{c|}{FETI} \\
		\(K\)   &      \(h\) &   \(\delta\) &   its &     time &its &    time &  \(L^{2}\) error &     RoC \\
		\hline
		6x6          & 0.004  &   0.008 &   343 & 0.052 & 40 & 0.103 & 3.47e-06 &         \\
		12x12        & 0.002  &   0.008 &   402 & 0.19  & 40 & 0.24 & 7.04e-07 & 2.3   \\
		24x24        & 0.001  &   0.008 &   398 & 0.89  & 40 & 0.76 & 1.77e-07 & 1.99 \\
		48x48        & 0.0005 &   0.008 &   405 & 4.79   & 41 & 1.98  & 4.45e-08 & 1.99 \\
		\hline
		\\
		\multicolumn{9}{c}{Fractional kernel \(\gamma_{s}\), \(s=0.4\)} \\
		\hline
		&&&\multicolumn{2}{c|}{CG}&\multicolumn{2}{c|}{FETI} \\
		\(K\)   &      \(h\) &   \(\delta\) &   its &     time &its &    time &  \(L^{2}\) error &     RoC \\
		\hline
		6x6          & 0.004  &   0.008 &   527 &  0.074 & 39 & 0.097 & 0.0014  &         \\
		12x12        & 0.002  &   0.008 &   705 &  0.28 & 36 & 0.21  & 7.26e-05 & 4.33 \\
		24x24        & 0.001  &   0.008 &   958 &  1.97  & 34 & 0.47  & 1.64e-05 & 2.14 \\
		48x48        & 0.0005 &   0.008 &  1197 & 11   & 32 & 1.59   & 3.85e-06 & 2.09 \\
		\hline
		\\
		\multicolumn{9}{c}{Peridynamic kernel \(\gammab_{PD}\)} \\
		\hline
		&&&\multicolumn{2}{c|}{CG}&\multicolumn{2}{c|}{FETI} \\
		\(K\)   &      \(h\) &   \(\delta\) &   its &     time &its &    time &  \(L^{2}\) error &     RoC \\
		\hline
		6x6          & 0.004  &   0.008 &   711 &  0.66 & 80 & 0.71 & 0.0089  &         \\
		12x12        & 0.002  &   0.008 &   824 &  2.17  & 73 & 1.56  & 9.34e-05 & 6.58 \\
		24x24        & 0.001  &   0.008 &   859 &  7.85  & 61 & 3.30  & 2.59e-05 & 1.85 \\
		48x48        & 0.0005 &   0.008 &   864 & 19.6   & 58 & 7.97  & 6.50e-06 & 1.99 \\
		\hline
	\end{tabular}
	
	\caption{Fixed horizon \(\delta=0.008\).
		From top to bottom: Constant kernel \(\gamma_{c}\), fractional kernel \(\gamma_{s}\) with \(s=0.4\) and peridynamic kernel \(\gammab_{PD}\).
	}
	\label{tab:fixedHorizon}
\end{table}

\subsection{Weak scaling for fixed \(\delta/h\)}
\label{sec:weak-scaling-ratio}

In a second set of experiments, we fix the ratio \(\delta/h = 4\).
Compared to the previous experiments, this means that the density of the discretization does not increase as the problem is refined.
We show the resulting iteration counts and timings in Table~\ref{tab:fixedRatio}.
It can be observed that the conjugate gradient solver does not scale for any of the considered kernels.
The FETI approach, on the other hand, converges in essentially constant number of iterations and almost perfect weak scaling.
We observe speedups of up to 7.9x (constant kernel), 14.6x (fractional kernel) and 8.9x (peridynamic kernel).
Based on our observations, we expect even larger improvements for more refined problems.

\begin{table}
	\centering
	
	\begin{tabular}{lrr|rr|rr}
		\multicolumn{7}{c}{Constant kernel \(\gamma_{c}\)} \\
		\hline
		&&&\multicolumn{2}{c|}{CG}&\multicolumn{2}{c}{FETI} \\
		\(K\)   &      \(h\) &   \(\delta\) &   its &     time &its &    time \\
		\hline
		6x6          & 0.002   &   0.008 &   402 &  1.19 & 39 & 1.12 \\ 
		12x12        & 0.001   &   0.004 &   733 &  2.25 & 39 & 1.19 \\ 
		24x24        & 0.0005  &   0.002 &  1436 &  5.44 & 40 & 1.27 \\ 
		48x48        & 0.00025 &   0.001 &  2799 & 12.7  & 40 & 1.61 \\ 
		\hline
		\\
		\multicolumn{7}{c}{Fractional kernel \(\gamma_{s}\), \(s=0.4\)} \\
		\hline
		&&&\multicolumn{2}{c|}{CG}&\multicolumn{2}{c}{FETI} \\
		\(K\)   &      \(h\) &   \(\delta\) &   its &     time &its &    time \\
		\hline
		6x6          & 0.002   &   0.008 &   706 &  1.91 & 40 & 1.11 \\ 
		12x12        & 0.001   &   0.004 &  1387 &  4.04 & 40 & 1.29 \\ 
		24x24        & 0.0005  &   0.002 &  2714 &  8.92 & 40 & 1.21 \\ 
		48x48        & 0.00025 &   0.001 &  5324 & 21.7  & 40 & 1.49 \\ 
		\hline
		\\
		\multicolumn{7}{c}{Peridynamic kernel \(\gammab_{PD}\)} \\
		\hline
		&&&\multicolumn{2}{c|}{CG}&\multicolumn{2}{c}{FETI} \\
		\(K\)   &      \(h\) &   \(\delta\) &   its &     time &its &    time \\
		\hline
		6x6          & 0.002   &   0.008 &   824 &  8.39 & 82 & 8.33 \\ 
		12x12        & 0.001   &   0.004 &  1582 & 18.12  & 84 & 9.43 \\ 
		24x24        & 0.0005  &   0.002 &  3093 & 41.87  & 87 & 9.72 \\ 
		48x48        & 0.00025 &   0.001 &  6063 & 87.85  & 87 & 9.84 \\ 
		\hline
	\end{tabular}
	
	\caption{
		Fixed ratio \(\delta/h=4\) of horizon and mesh size.
		From top to bottom: Constant kernel \(\gamma_{c}\), fractional kernel \(\gamma_{s}\) with \(s=0.4\) and peridynamic kernel \(\gammab_{PD}\).
		Optimal strong scaling as dashed lines.
	}
	\label{tab:fixedRatio}
\end{table}

\subsection{Strong scaling}
\label{sec:strong-scaling}

\begin{figure}
	\centering
	\includegraphics{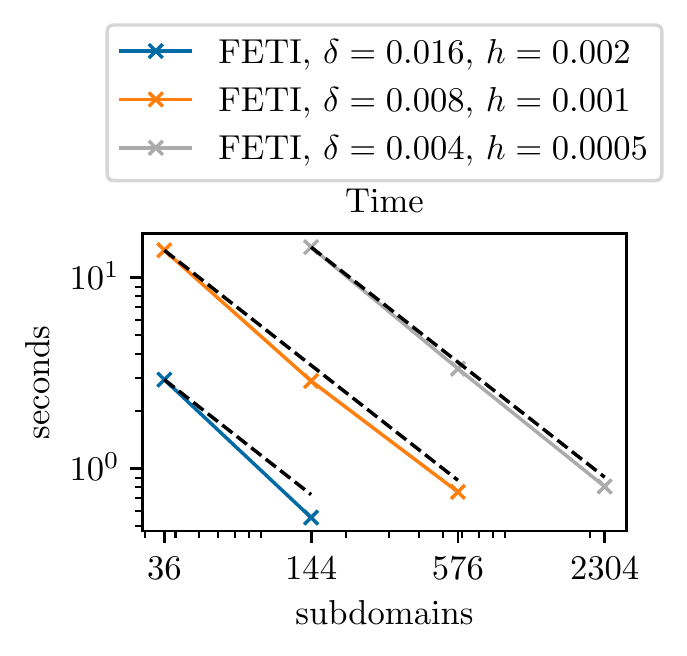}
	\includegraphics{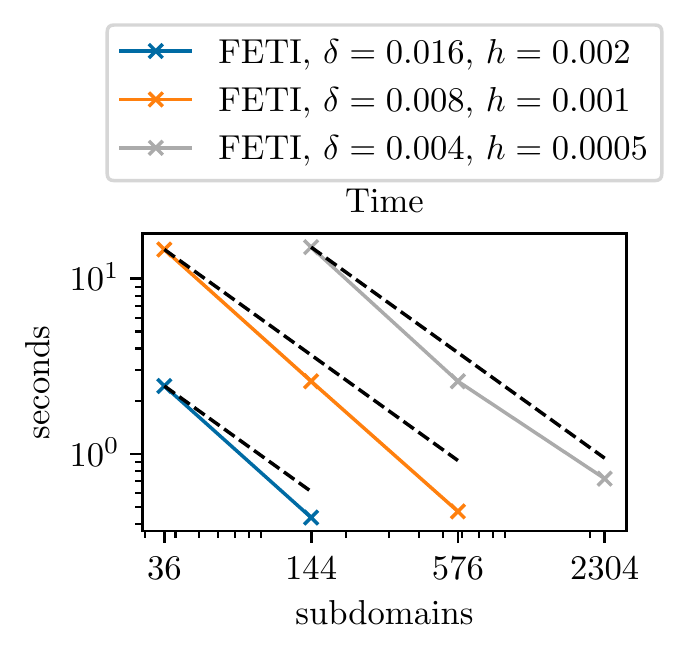}
	\includegraphics{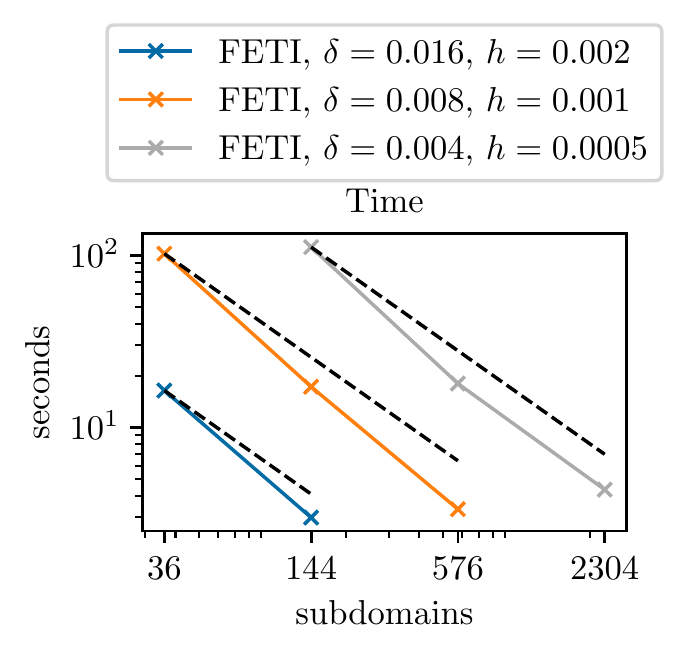}
	\caption{
		Strong scaling of FETI for \(\delta/h=8\).
		From top to bottom: Constant kernel \(\gamma_{c}\), fractional kernel \(\gamma_{s}\) with \(s=0.4\) and peridynamic kernel \(\gammab_{PD}\).
	}
	\label{fig:strongScaling}
\end{figure}

In a final set of experiments, we explore the strong scalability of the proposed FETI solver.
We fix \(\delta/h=8\), take \(h\in\{0.002,0.001,0.0005\}\) and solve the three problems on an increasing number of subdomains.
We display the solve timings for the three considered kernels in Figure~\ref{fig:strongScaling}, along with dashed lines for ideal strong scaling, i.e. \(\text{solve time}\sim 1/K\).
It can be observed that the timings scale even better than the ideal scaling.
This is explained by the fact that the subdomain problems are solved using direct solvers (CHOLMOD) which does not scale linearly in the size of the subdomain matrix.

\section{Conclusion}
\label{sec:conclusion}

The availability of efficient solvers is one of the main concerns for the use of nonlocal models in real world engineering applications.
In the present work we have developed a domain decomposition solver for finite element discretizations of nonlocal, finite horizon problems of integrable and fractional type.
While inspired by substructuring methods for PDEs which use nonoverlapping mesh partitions, as a consequence of nonlocality our approach requires overlapping subdomains.
We have covered the efficient implementation of the method using openly available scientific math libraries and demonstrated its weak and strong scalability with a series of numerical experiments.

\section{List of symbols}
\begin{table}[ht!]
		\begin{tabular}{lll}
                $B_{\delta, p}(\xb)$    & \pageref{inline:Bdeltap} & $\ell^p$-ball around $\xb$ \\
			$\widehat{\Omega}$														&					p. \pageref{inline:widehatOmega}				& given open, bounded domain \\
			$\interactionDomain$											&					\eqref{interdomain}				& $\delta$ strip around $\widehat{\Omega}$ \\
			$\Gamma^D$												&			p. \pageref{inline:nonlocalBoundary}						& nonlocal Dirichlet boundary \\
			$\Gamma^N$												&			p. \pageref{inline:nonlocalBoundary}							& nonlocal Neumann boundary (possibly empty)\\
			${\Omega}$														&					p. \pageref{inline:Omega}				& $\widehat{\Omega} \cup \Gamma^N$\\
                $\Gamma_h^D, \Gamma_h^N, \Omega_h$ & p. \pageref{inline:GammaDh}& geometric approximations of $\Gamma^D$, $\Gamma^N$ and $\Omega$\\
			$\mathcal{T}, \mcT^{\widehat\Omega}, \mcT^D, \dots$				&		\eqref{gridsgrids}							&triangulation of the domain\\
            $B_{\delta, 2}^h(\xb)$ & \eqref{eq:scalarkernel-h} & geometric approximation of the $\ell^2$-ball \\
			$\widehat{\Omega}_{k}$										&	Def \ref{def:nonlocalsubdivision}								& subdomain $\widehat{\Omega}_{k} \subset \Omega$ \\
			$\Gamma^D_{k}$										&	Def \ref{def:nonlocalsubdivision}									& subdomain $\Gamma^D_{k} \subset \Gamma^D$ \\
			$\Gamma_{k \ell}$								&\eqref{eq:subdomainDefinitions} 	& mutual interface $\widehat{\Omega}_k \cap \widehat{\Omega}_\ell$ for $k\neq \ell$\\
			$\Gamma_k$											&\eqref{eq:subdomainDefinitions}	& $\bigcup_{k\neq \ell} \Gamma_{k \ell}$, interface of $\widehat{\Omega}_k$
			\\
			$\Gamma^D_k$										&\eqref{eq:subdomainDefinitions}	& Dirichlet boundary of $\widehat{\Omega}_{k}$
			(possibly empty) \\
			$\Omega_k$													&\eqref{eq:subdomainDefinitions}	& interior subdomain $\widehat{\Omega}_k
			\setminus \Gamma_k$ \\
			$\Gamma$												&\eqref{eq:subdomainDefinitions} 	&  interface domain $\bigcup_{k=1}^K \Gamma_k$		\\
			$\Omega^h$														& Def \ref{def:femSpace}			& Nodes of a triangualation $\mcT$ of $\Omega$. \\
		\end{tabular}
	\caption{Sets and domains}
\end{table}

\begin{table}[ht!]
  \begin{tabular}{lll}
    $W(\Omega\cup\Gamma^{D})$				& Def \ref{def:energySpace} 	       & energy space on domain $\Omega\cup\Gamma^{D}$                                                                 \\
    $W^0(\Omega, \Gamma^D)$		& Def \ref{def:energySpace}	       & constrained energy space $w|_{\Gamma^D} = 0$                                                         \\
    \(W'(\Omega,\Gamma^{D})\)           & p. \pageref{inline:dualSpace}        & dual space of \(W^{0}(\Omega,\Gamma^{D})\) \\
    \(W^{t}(\Gamma^{D})\)           & p. \pageref{inline:traceSpace}        & trace space \\
    $W_h({\Omega_h\cup\Gamma^{D}_{h}})$	& Def \ref{def:femSpace}	       & finite element space                                                                                 \\
    $U_{\Omega, k}$			& p. \pageref{inline:U}		       & subspace of coefficients of $W^{0}_h({\widehat{\Omega}_k,\Gamma_{k}^{D}})$                                                       \\
    $U_{\Gamma, k}$			& p. \pageref{inline:U}		       & subspace of coefficients of $W^{0}_h({\widehat{\Omega}_k,\Gamma_{k}^{D}})$                                                       \\
    $U$					& p. \pageref{inline:U}		       & space of (possibly infeasible) solutions, $U = \prod_{k=1}^K U_{\Gamma, k}$                          \\
    $\Lambda = \R^{M_C}$		& \eqref{eq:constraintSpace}           & space of constraints, $M_C = \sum_{\xb_\ell \in \Gamma^h} \kerneldim (\zeta(\xb_\ell, \xb_\ell) -1)$ \\
    $\Lambda_{ad}$			& p. \pageref{inline:admissibleIncrements} & admissible increments of, $\Lambda_{ad} = \Pb(\Lambda)$
  \end{tabular}
  \caption{Spaces}
\end{table}

\begin{table}[ht!]
  \begin{tabular}{lll}
    $\mcA(u,v)$           & \eqref{eq:bil1}             & bilinear form for $u \in W({\widehat{\Omega}}\cup\Gamma^{D})$, $v \in W^0({\widehat{\Omega}, \Gamma^D})$                          \\
    $\mcF(v)$             & \eqref{eq:bil1}             & forcing term for $v \in W^0({\widehat{\Omega}, \Gamma^D})$                                                                        \\
    $\mcE$                & \eqref{singleenergy}        & energy functional                                                                                                                 \\
    \(\gamma_{h}\)        & \eqref{eq:scalarkernel-h}   & kernel with approximate interaction ball                                                                                          \\
    \(\mcA_{h}(u,v)\)     & p. \pageref{inline:bilQuad} & evaluation of \(\mcA\) using the kernel approximation \(\gamma_{h}\) and numerical quadrature                                     \\
    $\zeta(\xb, \yb)$     & \eqref{eq:zeta}   & counts overlaps of pairs                                                                                                          \\
    $\mcA_k(u_k,v_k)$     & \eqref{eq:bilinear-multi1}  & subdomain bilinear form for $u_k \in W({\widehat{\Omega}_{k}\cup\Gamma^D_{k}})$,  $v_k \in W^0({\widehat{\Omega}}_k, \Gamma^D_k)$ \\
    $\mcA_{k,h}(u_k,v_k)$ & \eqref{eq:bilinear-multi1}  & subdomain bilinear form using the kernel approximation \(\gamma_{h}\) and numerical quadrature                                    \\
    $\mcF_k(v_k)$         & \eqref{eq:bilinear-multi5}  & subdomain forcing term for $v_k \in W^0({\widehat{\Omega}}_k, \Gamma^D_k)$                                                        \\
    $\mcE_k$              & \eqref{eq:subE}             & subdomain energy functional
  \end{tabular}
\caption{Functions}
\end{table}
\begin{table}[ht!]
	\begin{tabular}{lll}
          $\mbB\colon U \rightarrow \Lambda$														&		\eqref{eq:constraintSpace}							& Constraints, $\ub \in$ ker$(\mbB)$ is feasible\\
          $\mbA^k$
                                                                                                                                                        &		\eqref{eq:subStiffnessMatrix}							&subdomain stiffness matrix \\
          $\mbA^k_{\Omega\Omega}$ 	& \eqref{eq:subStiffnessMatrixSubBlocks}		&submatrix in $\Omega_k$ \\
		$\mbA^k_{\Gamma\Gamma}$ 		&	\eqref{eq:subStiffnessMatrixSubBlocks}								&submatrix in $\Gamma_k$ \\
		$\mbS_k$					& \eqref{eq:subSchurComplement}		& subdomain Schur complement\\
		$\mbS \colon U \rightarrow U$															& \eqref{eq:SchurComplement}		& $\diag(\mbS_k)$ \\
		$\ub^h \in U$																									&	\eqref{eq:SchurComplement}								& coefficients $(\ub^{h,1}_{\Gamma}, \dots, \ub^{h,K}_{\Gamma})$ \\
		$\fb^h \in U$	& 		\eqref{eq:SchurComplement}
		& reduced forcing term  \\
		$\Zb: \R^{\NFloating} \rightarrow U$															& 		\eqref{eq:matKerS}							& Null space of $\mbS$ \\
		$\Gb:  \R^{\NFloating} \rightarrow \Lambda$														&		p. \pageref{inline:G}							& $\Bb \Zb$ \\
		$\Mb^{-1} = (\Bb  \Db^{-1} \Bb ^\top )^{-1}  \Bb \Db^{-1} ~ \mbS ~ \Db^{-1} \Bb^\top  (\Bb \Db^{-1}  \Bb^\top )^{-1}$																						& \eqref{eq:preconditioner}			& Dirichlet preconditioner \\
		$\Db$																									& \eqref{eq:preconditioner}	& weights	\\
		$\Pb = \Ib - \Gb(\Gb^\top \Gb)^{-1} \Gb^\top : \Lambda \rightarrow \Lambda_{ad}$			&	\eqref{eq:P}								& orthogonal projection \\
		$\alpha \in \R^{\NFloating}$																			&	\eqref{eq:forAlpha}								& dofs due to rigid body modes\\
		$\Fb$																									&	\eqref{eq:FETIproblem}		& system $\Bb \mbS^+ \Bb^\top$ \\
		$\db$								&	\eqref{eq:FETIproblem}								& $\Bb \mbS^+ \fb$ \\
		$\eb$	&\eqref{eq:FETIproblem}		& $\Zb^\top \fb$
	\end{tabular}
	\caption{Matrices and vectors}
\end{table}

\section{Funding}
\label{sec:funding}

This work was supported by the Laboratory Directed Research and Development program at
Sandia National Laboratories, a multimission laboratory managed and operated by National
Technology and Engineering Solutions of Sandia LLC, a wholly owned subsidiary of Honeywell
International Inc. for the U.S. Department of Energy’s National Nuclear Security Administration
under contract DE-NA0003525.

\bibliography{domdecomp.bib}

\end{document}